%% file: main.tex
\documentclass{article}
\usepackage[utf8]{inputenc}
\usepackage[round, sort]{natbib}

\usepackage{amsmath}
\usepackage{amssymb}
\usepackage{amsthm}
\usepackage{mathtools}
\usepackage[mathscr]{euscript}

\usepackage{pgf,tikz,pgfplots}

\usepackage{bbm}

\usepackage{tikz-cd}
\usetikzlibrary{patterns}

\usetikzlibrary{arrows,decorations.markings}
\usetikzlibrary{math}

\usepackage{dirtytalk}

\usepackage{enumerate}
\usepackage[shortlabels]{enumitem}

\usepackage{subfigure}

\usepackage[margin=1in]{geometry}

\usepackage{xcolor}

\usepackage{setspace}
\theoremstyle{definition}
\newtheorem{definition}{Definition}

\theoremstyle{plain}
\newtheorem{thm}{Theorem}

\theoremstyle{plain}
\newtheorem{lemma}{Lemma}

\theoremstyle{plain}
\newtheorem{prop}{Proposition}

\theoremstyle{plain}
\newtheorem{cor}{Corollary}

\usepackage[T2A,T1]{fontenc}

\usepackage[symbol]{footmisc}

\usepackage{array}
\usepackage{makecell}

\renewcommand{\d}{\,\mathrm{d}}
\newcommand{\limn}{\lim_{n\rightarrow \infty}}

\newcommand{\tor}{\mathbb{T}^2}

\usepackage{subfiles}

\usepackage{multicol}
\providecommand{\keywords}[1]
{
  \small	
  \textbf{\textit{Keywords---}} #1
}

\providecommand{\Acknowledgements}[1]
{
  \small	
  \textbf{\textit{Acknowledgements---}} #1
}

\title{Exponential Mixing by Orthogonal Non-Monotonic Shears}
\author{J. Myers Hill$^{1\,\star}$, R. Sturman$^{2}$, M. C. T. Wilson$^{3}$  \\
        \small $^{1}$EPSRC CDT in Fluid Dynamics, University of Leeds, Leeds LS2 9JT, United Kingdom \\
        \small $^{2}$School of Mathematics, University of Leeds, Leeds LS2 9JT, United Kingdom \\
        \small $^{3}$School of Mechanical Engineering, University of Leeds, Leeds LS2 9JT, United Kingdom \\
        \small $^\star$ E: scjmh@leeds.ac.uk
}

\pgfplotsset{compat=1.17}
\begin{document}

\date{}

\maketitle

\begin{abstract}
   Non-monotonic velocity profiles are an inherent feature of mixing flows obeying non-slip boundary conditions. There are, however, few known models of laminar mixing which incorporate this feature and have proven mixing properties. Here we present such a model, alternating between two non-monotonic shear flows which act in orthogonal (i.e. perpendicular) directions. Each shear is defined by an independent variable, giving a two-dimensional parameter space within which we prove the mixing property over open subsets. Within these mixing windows, we use results from the billiards literature to establish exponential mixing rates. Outside of these windows, we find large parameter regions where elliptic islands persist, leading to poor mixing. Finally, we comment on the challenges of extending these mixing windows and the potential for a non-exponential mixing rate at particular parameter values.
   
\end{abstract}


\Acknowledgements{JMH is supported by EPSRC under Grant Ref. EP/L01615X/1.}

\section{Introduction}
\label{sec:intro}
\subfile{sections/intro}

\section{Proof Outline}

\subfile{sections/outline}
\label{sec:outline}

\section{Establishing non-uniform hyperbolicity}

\label{sec:lyp}
\subfile{sections/lyp}

\section{Establishing ergodicity}
\label{sec:man}
\subfile{sections/man}

\section{Establishing the Bernoulli property}
\label{sec:repMan}
\subfile{sections/repMan}

\section{Rate of mixing}
\label{sec:rates}
\subfile{sections/mixingRate}

\section{The two dimensional parameter space}
\label{sec:param}
\subfile{sections/2Dparam}

\section{Special case}
\label{sec:otm}
\subfile{sections/OTM}

\section{Discussion}
\label{sec:discussion}
\subfile{sections/discussion}

\section*{Appendix}
\subfile{sections/appendix}

\newpage

\bibliography{refs.bib}
\bibliographystyle{apalike}

\end{document}

%% file: sections/intro.tex
\subsection{Background}

Mixing a fluid in some domain $X$ by chaotic advection, essentially stirring, typically concerns the study of time-$T$ periodic incompressible laminar flows $\mathbf{v}(\textbf{x},t) = \mathbf{v}(\textbf{x},t+T)$, $\textbf{x} \in X$. The dynamical features of these flows are described by the trajectories of fluid particles within the system, in particular their positions after each time period $T$. This defines a map $f:X \rightarrow X$ which sends the initial position of a particle to its position after time $T$. The long term behaviour of the flow $\mathbf{v}$ on $X$ is then described by repeated iterations of the map $f$, denoted by $f^n$, and the incompressibility condition on $\mathbf{v}$ tells us that $f$ preserves the Lebesgue measure $\mu$ on $X$. Ergodic theory provides the mathematical framework for understanding the long term behaviour of iterating such maps, most importantly it gives a precise definition for what it means for a map to be mixing:
\begin{definition}
A $\mu$-preserving invertible map $f:X \rightarrow X$ is \emph{mixing} if $\limn \mu(f^{n}(B) \cap A) = \mu(A)\mu(B)$ for all measurable $A,B\subset X$, $\mu(X)=1$.
\end{definition}
A hallmark of chaotic advection in two dimensions is an iterative `stretching and folding' type action. Starting with a blob of fluid $B$, this action transforms $f^n(B)$ into a thin fluid filament which is repeatedly folded to spread across the entire domain, or any target set $A$. A natural model for this behaviour is to compose shear maps, which stretch while preserving $\mu$, on the torus $\tor$, whose periodic boundaries interweave the long fluid filaments. A canonical example is the Cat Map $H:\tor \rightarrow \tor$ (\citealp{arnold_ergodic_1968}). Parameterise $\tor$ by $(x,y) \in (\mathbb{R}/\mathbb{Z})^2$, then $H = G \circ F$ composes horizontal and vertical shears, written in matrix form as
\[ \begin{pmatrix} x \\ y \end{pmatrix} \mapsto \underbrace{\begin{pmatrix} 1 & 0 \\ 1 & 1 \end{pmatrix}}_{DG} \underbrace{\begin{pmatrix} 1 & 1 \\ 0 & 1 \end{pmatrix}}_{DF}  \begin{pmatrix} x \\ y \end{pmatrix} \text{ mod 1} = \underbrace{\begin{pmatrix} 1 & 1 \\ 1 & 2 \end{pmatrix}}_{M} \begin{pmatrix} x \\ y \end{pmatrix} \text{ mod 1}, \]
where $DF$, $DG$ denote the Jacobians of the maps $F$, $G$. Since $H$ can be defined using a single hyperbolic matrix $M$, the map is \emph{uniformly hyperbolic}, with the same magnitude and directions of expansion, contraction across the entire domain. While this allows for a straightforward proof of the mixing property, mixing in a realistic fluid flow is typically non-uniform due to the influence of walls and non-linear velocity profiles. The key barriers to mixing in this setting are \emph{elliptic islands}, invariant subsets within which particle paths trace out closed curves. From a fluids perspective these curves form material lines in the flow which particles cannot penetrate (except by diffusion), leading to poor mixing (\citealp{ottino_kinematics_1989}). An illustration of an island pair is given later in Figure \ref{fig:islandsPoin}. Previous studies which try to minimise island structures in realistic mixing flows include \cite{franjione_symmetry_1992}, looking at eggbeater and duct flows, and \cite{hertzsch_dna_2007}, looking at mixing in DNA microarrays. 

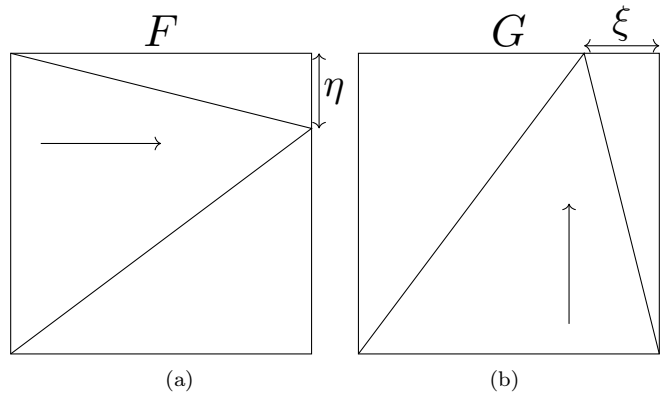
\begin{figure}[h]
    \centering
\subfigure[][]{%
\label{fig:Feta}%
    \begin{tikzpicture}
    \draw (0,0) rectangle (4,4);
\draw (0,0)--(4,3)--(0,4);
\draw[<->] (4.1,4)--(4.1,3);
\node[scale=1.5] at (4.3,3.5) {$\eta$};
\draw[->] (0.4,2.8)--(2,2.8);
\node[scale=1.7] at (2,4.3) {$F$};
    \end{tikzpicture}}%
\subfigure[][]{%
\label{fig:Gxi}%
    \begin{tikzpicture}
    \draw (0,0) rectangle (4,4);
\draw (0,0)--(3,4) -- (4,0);
\node[scale=1.5] at (3.5,4.4) {$\xi$};
\draw[<->] (3,4.1)--(4,4.1);
\draw[->] (2.8,0.4)--(2.8,2);
\node[scale=1.7] at (2,4.3) {$G$};
    \end{tikzpicture}}%

    \caption{A family of area preserving maps $H_{(\xi,\eta)} = G_\xi \circ F_\eta$ parameterised by $0<\eta,\xi<1$.}
    \label{fig:mapDefn}
\end{figure}

Several mappings exist which incorporate realistic flow phenomena and still allow for a proof of the mixing property. Linked Twist Maps (\citealp{burton_ergodicity_1980}), hereafter LTMs, compose monotonic shears which act on annuli of the torus, leaving a region invariant which models a boundary within the domain. Mixing properties can be shown in the co-rotating case, where the shears act in the positive $x$ and $y$ directions (\citealp{wojtkowski_linked_1980}). In the counter-rotating case, where the direction of one shear is reversed, island structures develop and can only be broken up by taking strong shears. This highlights a potential challenge for mixing by non-monotonic shears, which inherently exhibit this counter rotating quality. Indeed, there are few examples of non-monotonic toral maps with proven mixing properties. \citeauthor{cerbelli_continuous_2005}'s map (hereafter the CG Map), studied in \cite{cerbelli_continuous_2005}, \cite{mackay_cerbelli_2006}, and \cite{cerbelli_characterization_2008}, incorporates non-monotonicity into the first shear by taking
\[ F(x,y) = \begin{cases} (x + 2y,y)\text{ mod 1}  & \text{for } y \leq \frac{1}{2}, \\
(x+2(1-y),y)\text{ mod 1} & \text{for } y \geq \frac{1}{2},
\end{cases}   \]
and leaves $G$ unchanged. While this introduces a non-hyperbolic derivative matrix over half the domain, a unique geometric feature of the map ensures that this does not compromise long term stretching behaviour. Indeed, hyperbolic and mixing properties of the CG Map can be proven by quite direct means, with analysis of the unstable foliation revealing that fluid filaments get stretched and folded in a very regimented fashion. This is not the case for generic non-monotonic shears $F$, as shown in \cite{myers_hill_continuous_2021}, where shears of the form
\begin{equation}
\label{eq:Feta}
    F_\eta(x,y) = \begin{cases}
\left(  x + \frac{1}{1-\eta} y , y  \right) \text{ mod 1 } & \text{ for } y \leq 1-\eta, \\

\left(  x + \frac{1}{\eta} (1-y) , y  \right) \text{ mod 1 } & \text{ for } y \geq 1-\eta, \\
\end{cases}
\end{equation}
are considered with $G$ left unchanged. An illustration of this shear is given in Figure \ref{fig:Feta}; note that parameter values $\eta=0$ and $\eta=1/2$ give the Cat Map and CG Map respectively. Mixing properties over subsets of the parameter space $0<\eta<1/2$ are shown using a scheme from \cite{katok_invariant_1986}, Theorem \ref{thm:katok-strelcyn} in the present work, which gives comparatively easy to verify conditions under which non-uniformly hyperbolic systems\footnote{In particular non-uniformly hyperbolic systems with singularities, which are inherent to piecewise linear maps.} are mixing, compared to arguing by direct means. Non-uniform hyperbolicity ensures the existence (almost everywhere) of local stable and unstable manifolds which, roughly speaking, describe the characteristic local flow direction in backwards and forwards time respectively (see for example \citealp{beigie_invariant_1994}). The way in which blobs of fluid are stretched and spread across the domain is then described well by the images of these local manifolds, and mixing properties follow from intersection conditions on these images. Broadly speaking, the key challenge to showing these conditions for non-monotonic systems is the sign alternating property as described in \cite{cerbelli_continuous_2005}. While hyperbolicity ensures that the images of local manifolds grow exponentially in length, non-linear shears like $F_\eta$ fold these images back on themselves, potentially inhibiting their spread across the domain (also commented on in \citealp{przytycki_ergodicity_1983}). This challenge can often be overcome by establishing sufficiently strong stretching behaviour over one or several iterates, a method we will employ here to prove mixing properties of maps composing two non-monotonic shears. While taking $G$ similar to $F_\eta$ (see Figure \ref{fig:Gxi}) means that the images of local manifolds will fold back on themselves more often, in turn the gradients defining $G$ will become steeper, giving stronger stretching behaviour.

Alongside knowing whether a map is mixing, it is desirable to know the rate at which we approach a mixed state. That is, what is the rate of decay of $\left|\mu(f^{n}(B) \cap A) - \mu(A)\mu(B)\right|$ with $n$? Taking indicator functions $\mathbbm{1}_A(x) = 1$ if $x \in A$, 0 otherwise, and defining the correlation function
\begin{equation}
    \label{eq:corrFunction}
    C_n(\phi,\psi,f,\mu) = \int \left(\phi \circ f^n \right) \psi \d \mu - \int \phi \d \mu \int \psi \d \mu 
\end{equation}
for observables $\phi,\psi:\tor \rightarrow \mathbb{R}$, this amounts to studying the decay of $|C_n(\mathbbm{1}_B,\mathbbm{1}_A,f^{-1},\mu)|$. For uniformly hyperbolic maps, proving exponential correlation decay rate is quite straightforward using a \emph{Young tower} construction, shown explicitly for the Cat Map in \cite{chernov_decay_2000}, based on the theory developed in \cite{young_statistical_1998}, \cite{young_recurrence_1999}. The theory has been extended to hyperbolic systems with singularities (\citealp{chernov_decay_1999}) culminating in schemes from \cite{chernov_billiards_2005}, which give conditions under which systems enjoy exponential or polynomial decay of correlations, provided hyperbolicity is sufficiently strong. While aimed at application to billiards systems, the scheme is readily applicable to models of fluid mixing, most recently in \cite{springham_polynomial_2014} where the mixing rate for a wide class of linked twist maps is shown to be at worst polynomial. 

\subsection{Statement of results}

Let $H_{(\xi,\eta)}:\tor \rightarrow \tor$ be the composition of two shears $F_\eta$ and $G_\xi$, where $F_\eta$ is given by (\ref{eq:Feta}) and $G_\xi$ maps
\[ (x,y) \mapsto 
\begin{cases}
\left(  x  , y + \frac{1}{1-\xi} x  \right) \text{ mod 1 } & \text{ for } x \leq 1-\xi, \\

\left(  x  , y + \frac{1}{\xi} (1-x) \right) \text{ mod 1 } & \text{ for } x \geq 1-\xi. \\
\end{cases}\]
\begin{figure}
    \centering
     \begin{tikzpicture}
     \tikzmath{\e = 0.2;} 
    \node[scale=1.2] at (0,-0.22) {
    \begin{tikzpicture}
    
    \fill[gray!40] (0,4-4*\e) rectangle (4-4*\e,4);
    
    \fill[gray!20] (4-4*\e,0) rectangle (4,4-4*\e);
    
    \draw (0,0) rectangle (4,4);
    \draw (4-4*\e,0) -- (4-4*\e,4);
    \draw (0,4-4*\e) -- (4,4-4*\e);

    \node[scale=0.8] at (4-4*\e,-0.2) {$1-\eta$};
    
    \node[scale=1.5] at (2-0.5*4*\e,2-0.5*4*\e) {$R_1$};
    \node[scale=1.5] at (4-0.5*4*\e,2-0.5*4*\e) {$R_2$};
    \node[scale=1.5] at (2-0.5*4*\e,4-0.5*4*\e) {$R_3$};
    \node[scale=1.5] at (4-0.5*4*\e,4-0.5*4*\e) {$R_4$};
    
    \end{tikzpicture}
    };
    \node[scale=1.2] at (-5.5,0.025) {
    \begin{tikzpicture}
    
    \filldraw[fill=gray!40] (0, {4*(1-\e+\e*\e)}) -- ({4*(1-\e)},4) -- (0,4) --  (0, {4*(1-\e+\e*\e)});
    
    \filldraw[fill=gray!40] (4,4) -- (0,{4*(1-\e)}) -- ({4*(1-\e)},{4*(1-\e)}) -- (4, {4*(1-\e+\e*\e)}) -- (4,4);
    
    \filldraw[fill=gray!20] ({4*(1-\e)},{4*(1-\e)}) -- (4,{4*(1-\e)}) -- (4, {4*(1-\e)^2}) -- ({4*(1-\e)},{4*(1-\e)});
    
    \filldraw[fill=gray!20] (0,{4*(1-\e)}) -- (4,0) -- ({4*(1-\e)},0) -- (0, {4*(1-\e)^2}) -- (0,{4*(1-\e)});
    
    \draw (0,0) rectangle (4,4);

    \node[scale=0.8] at (-0.4,4-4*\e) {$1-\eta$};
    
        \node[scale=0.8] at (4.2, {4*(1-\e+\e*\e)}) {$y_2$};
        \node[scale=0.8] at (4.2, {4*(1-\e)^2}) {$y_1$};
    \node at (2-0.5*4*\e,2-0.5*4*\e) {$A_2$};
    \node at (4-0.3*4*\e, {4*(1-1.25*\e)} ) {$A_2$};
    
    \node at (4-0.5*4*\e,2-0.5*4*\e) {$A_1$};

    \node at (4-0.5*4*\e,4-0.5*4*\e) {$A_3$};
    \node at (0.5*4*\e,4-0.3*4*\e) {$A_3$};
    \node at (2*\e,2*\e) {$A_1$};
    \end{tikzpicture}
    };
    
    \node[scale=1.2] at (5,-0.2) {
    \begin{tikzpicture}
    
    \filldraw[fill = gray!40] (0,4-4*\e) -- ({4*\e*(1-\e)},4) -- (0,4) -- (0,4-4*\e);
    \filldraw[fill = gray!40] (0,0) -- (4-4*\e,4) -- (4-4*\e,4-4*\e) -- ({4*\e*(1-\e)},0) -- (0,0);
    
    \filldraw[fill=gray!20] (4-4*\e,0) -- (4-4*\e,4-4*\e) -- (4-4*\e*\e,0) -- (4-4*\e,0);
    \filldraw[fill=gray!20] (4,0) -- (4-4*\e,4) -- (4-4*\e*\e,4)-- (4,4-4*\e) -- (4,0);
    
    \draw (0,0) rectangle (4,4);
    
    \node[scale=0.8] at ({4*\e*(1-\e)},-0.2) {$x_1$};
    \node[scale=0.8] at ({4*(1-\e*\e)},-0.2) {$x_2$};
    

    \node at ({0.6*4*(1-\e)},2) {$A_3'$};
    \node at ({0.7*4*(1-\e)},1) {$A_1'$};
    \node at ({0.4*4*(1-\e)},3) {$A_1'$};  
        
    \node at (0.25*4*\e,4-0.3*4*\e) {$A_3'$};
     \node at (4-0.3*4*\e, {4*(1-1.25*\e)} ) {$A_2'$};
    \node at ({4-2.5*\e}, 2*\e ) {$A_2'$};
    \end{tikzpicture}
    };

    \end{tikzpicture}
    
    \caption{A partition of the torus into four rectangles $R_j$, and their images $A_j$ under $F_\eta^{-1}$, $A_j'$ under $G_\xi$. The smallest partition elements $A_4$ and $A_4'$ are left unlabelled. Case illustrated $\xi=\eta=0.2$.}
    \label{fig:Apartition}
\end{figure}
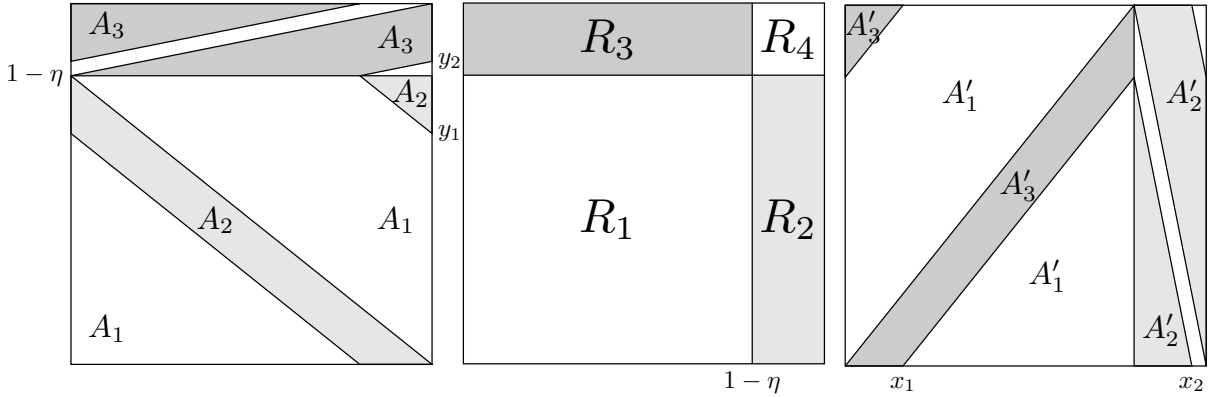
Write $H = H_{(\xi,\eta)}$ and partition the torus into four rectangles $R_j$ using the lines $x=0$, $y=0$, $x=1-\xi$, and $y=1-\eta$, as shown in Figure \ref{fig:Apartition}. Letting $A_j = F_\eta^{-1}(R_j)$, the derivative $DH$ is defined everywhere outside of the set $\mathcal{D} = \cup_j \partial A_j$ and is constant on each of the $A_j$. The matrices $M_j =  DH|_{A_j}$ are given by 
\[ M_1 = \begin{pmatrix} 1 & \frac{1}{1-\eta} \\ \frac{1}{1-\xi} & 1+\frac{1}{(1-\xi)(1-\eta)} \end{pmatrix}, \quad M_2 = \begin{pmatrix} 1 & \frac{1}{1-\eta} \\ -\frac{1}{\xi} & 1-\frac{1}{\xi(1-\eta)} \end{pmatrix}, \] 
\[ M_3 = \begin{pmatrix} 1 & -\frac{1}{\eta} \\ \frac{1}{1-\xi} & 1-\frac{1}{\eta(1-\xi)} \end{pmatrix}, \quad M_4 = \begin{pmatrix} 1 & -\frac{1}{\eta} \\ -\frac{1}{\xi} & 1+\frac{1}{\eta\xi} \end{pmatrix}. \]
Letting $A_j' = G(R_j) = H(A_j)$, the derivative of $H^{-1}$ is defined everywhere outside of the set $\mathcal{D}' = \cup_j \partial A_j'$ and is constant on each of the $A_j'$. The labelled intersections with the axes are $y_1 = (1-\xi)(1-\eta)$, $y_2 = 1-\eta+\xi\eta$, $x_1 = \eta(1-\xi)$, and $x_2 = 1-\xi\eta$.

\begin{figure}%
\centering
\subfigure[][]{%
\label{fig:expPoin}%
\includegraphics[width=0.33\linewidth]{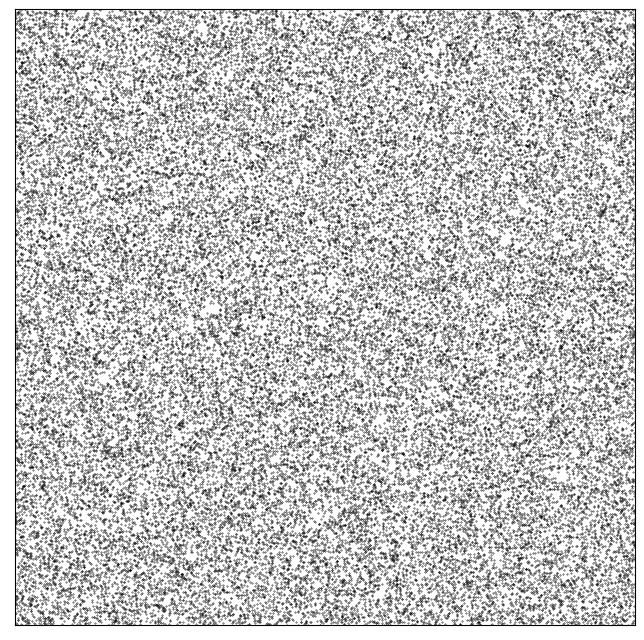}}%
\subfigure[][]{%
\label{fig:islandsPoin}%
\includegraphics[width=0.33\linewidth]{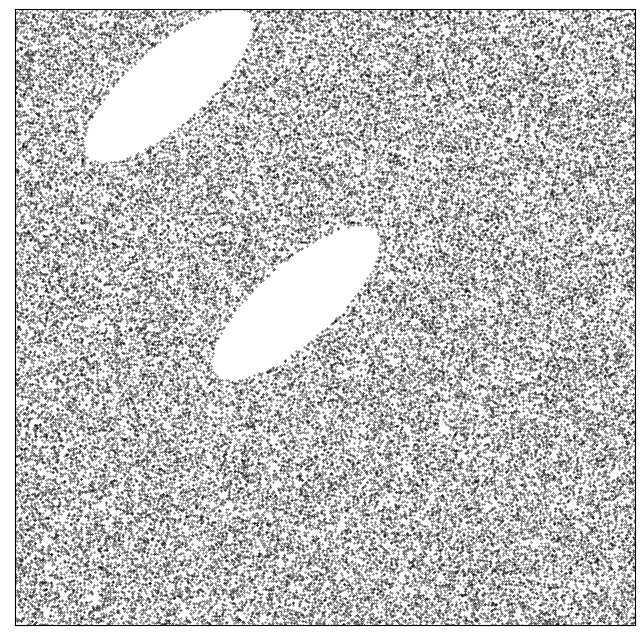}}%
\subfigure[][]{%
\label{fig:otmPoin}%
\includegraphics[width=0.33\linewidth]{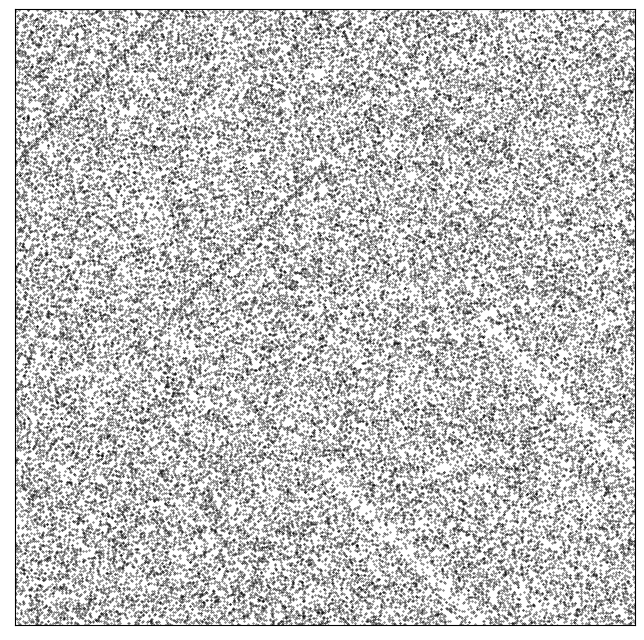}}%

\caption{Poincar\'e sections for $H_{(\xi,\eta)}$ at parameter values (a): $\xi=\eta=1/10$, (b): $\xi=3/10,\eta =6/10$, (c): $\xi=\eta=1/2$. Starting point is $z_0=(1/\sqrt{2},1/\sqrt{3})$, 50,000 iterates shown. In (a) we see fully ergodic behaviour, in (b) island structures are present, and in (c) orbits can become trapped near invariant sets of line segments.}%
\label{fig:poincare}%
\end{figure}

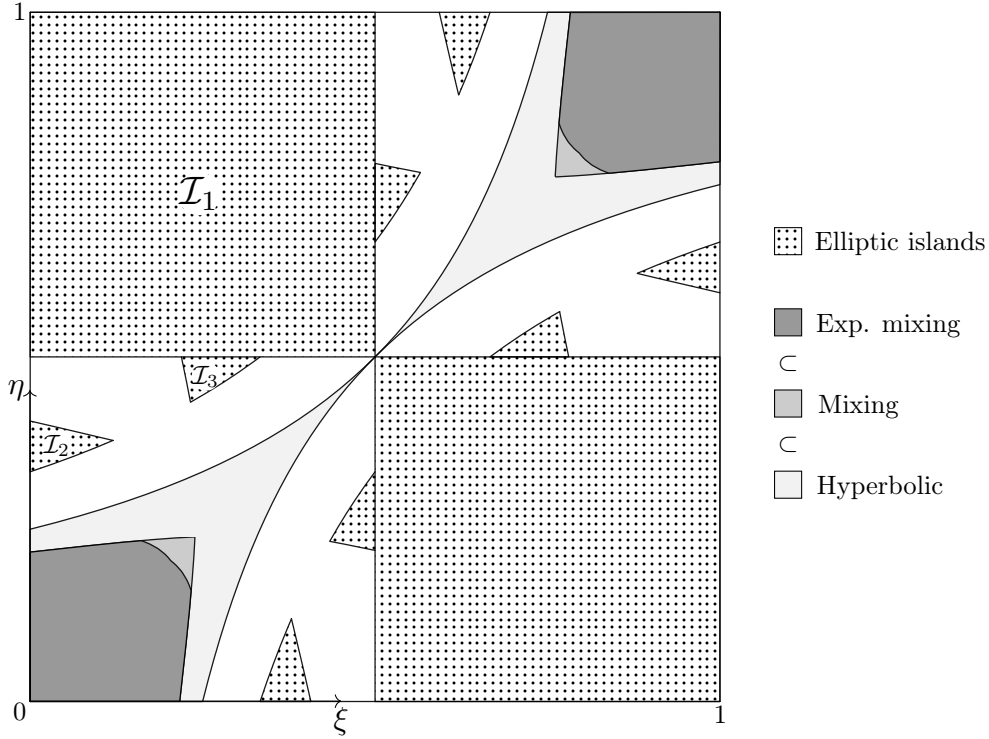
\begin{figure}
    \centering
    \begin{tikzpicture}[scale=1.8]
    \node at (2.487,2.542) {
    
    \definecolor{white}{RGB}{255,255,255}

\def \globalscale {0.231*18/25}
\begin{tikzpicture}[y=0.80pt, x=0.80pt, yscale=-\globalscale, xscale=\globalscale, inner sep=0pt, outer sep=0pt]
\begin{scope}[scale=2.000]
    \path[fill=white,rounded corners=0.0000cm] (0.0000,0.0000) rectangle (1000.0000,1000.0000);

        \path[draw=black,draw opacity=0.900,line cap=round,line join=round,line width=0.500pt,miter limit=10.00] (9.3039,989.6961) -- (9.3039,989.6961) -- (9.3039,9.3039);

        \path[draw=black,draw opacity=0.900,line cap=round,line join=round,line width=0.500pt,miter limit=10.00] (9.3039,9.3039) -- (9.3039,9.3039) -- (989.6961,9.3039);

        \path[draw=black,draw opacity=0.900,line cap=round,line join=round,line width=0.500pt,miter limit=10.00] (9.3039,989.6961) -- (9.3039,989.6961) -- (989.6961,989.6961);

        \path[draw=black,draw opacity=0.900,line cap=round,line join=round,line width=0.500pt,miter limit=10.00] (989.6961,989.6961) -- (989.6961,989.6961) -- (989.6961,9.3039);

        \path[fill=gray!10,draw=black,draw opacity=0.900,line cap=round,line join=round,line width=0.500pt,miter limit=10.00] (9.3039,989.6961) -- (9.3039,744.5980) -- (9.3039,744.5980) -- (45.3984,735.2295) -- (79.8223,725.6021) -- (112.5371,715.7524) -- (143.5430,705.7141) -- (173.3281,695.3531) -- (201.6484,684.7752) -- (228.7480,673.9160) -- (254.6270,662.7987) -- (279.5293,651.3358) -- (303.4551,639.5389) -- (326.6484,627.2905) -- (348.8652,614.7263) -- (370.3496,601.7191) -- (391.1016,588.2688) -- (411.3652,574.2035) -- (430.8965,559.6811) -- (449.6953,544.7112) -- (467.5176,529.5236) -- (484.3633,514.1833) -- (500.2324,498.7665) -- (515.3691,483.0999) -- (529.7734,467.2339) -- (543.6895,450.9323) -- (557.1172,434.2085) -- (570.0566,417.0802) -- (582.5078,399.5705) -- (594.4707,381.7083) -- (605.9453,363.5287) -- (617.1758,344.6515) -- (627.9180,325.4984) -- (638.4160,305.6487) -- (648.6699,285.0809) -- (658.6797,263.7737) -- (668.4453,241.7066) -- (677.9668,218.8599) -- (687.2441,195.2155) -- (696.2773,170.7566) -- (705.0664,145.4686) -- (713.6113,119.3394) -- (722.1562,91.5412) -- (730.4570,62.7825) -- (738.5137,33.0519) -- (744.5980,9.3039) -- (989.6961,9.3039);

        \path[fill=gray!10,draw=black,draw opacity=0.900,line cap=round,line join=round,line width=0.500pt,miter limit=10.00] (9.3039,989.6961) -- (254.4020,989.6961) -- (254.4020,989.6961) -- (262.4615,958.4844) -- (270.7092,928.5365) -- (279.1173,899.8906) -- (287.7688,872.2214) -- (296.5457,845.8542) -- (305.5369,820.4635) -- (314.7295,796.0495) -- (324.1076,772.6120) -- (333.7944,749.8255) -- (343.6413,728.0156) -- (353.7828,706.8568) -- (364.2171,686.3490) -- (374.9407,666.4922) -- (385.9478,647.2865) -- (397.2300,628.7318) -- (408.7761,610.8281) -- (420.8012,593.2500) -- (433.3289,575.9974) -- (446.1258,559.3958) -- (459.4445,543.1198) -- (473.3091,527.1693) -- (487.7446,511.5443) -- (502.7770,496.2448) -- (518.7846,480.9453) -- (535.4903,465.9714) -- (552.9281,451.3229) -- (571.5612,436.6745) -- (591.0581,422.3516) -- (611.4642,408.3542) -- (632.3005,395.0078) -- (654.0668,381.9870) -- (676.2019,369.6172) -- (699.2654,357.5729) -- (723.2900,345.8542) -- (748.3075,334.4609) -- (774.3481,323.3932) -- (801.4400,312.6510) -- (829.6085,302.2344) -- (858.8754,292.1432) -- (890.3078,282.0521) -- (923.0231,272.2865) -- (957.0435,262.8464) -- (989.6961,254.4020) -- (989.6961,9.3039);

        \path[fill=gray!40,draw=black,draw opacity=0.900,line cap=round,line join=round,line width=0.500pt,miter limit=10.00] (221.8336,989.4902) -- (221.8336,989.4902) -- (230.6817,908.9433) -- (235.8281,857.9290) -- (239.5646,815.9062) -- (241.9118,784.1282) -- (243.4491,756.5801) --     (242.4199,755.5509) -- (242.4199,755.5509) -- (214.8717,757.0882) -- (183.0938,759.4355) -- (141.0710,763.1719) -- (90.0567,768.3183) -- (9.5098,777.1664) -- (9.5098,989.4902);

        \path[fill=gray!40,draw=black,draw opacity=0.900,line cap=round,line join=round,line width=0.500pt,miter limit=10.00] (777.1664,9.5098) -- (777.1664,9.5098) -- (769.6814,77.6250) -- (764.2275,131.4150) -- (760.6707,170.5180) -- (757.4084,212.4865) -- (755.5602,242.4199) -- (756.5801,243.4491) -- (756.5801,243.4491) -- (784.1283,241.9118) -- (803.0731,240.5731) -- (803.0731,240.5731) -- (835.5474,237.8911) -- (886.2188,233.0369) -- (956.5313,225.5233) -- (989.4902,221.8336) -- (989.6961,9.3039);

        \path[fill=gray!80,draw=black,draw opacity=0.900,line cap=round,line join=round,line width=0.500pt,miter limit=10.00](9.5098,989.4902) -- (221.8336,989.4902) -- (221.8336,989.4902) -- (230.6817,908.9433) -- (235.8281,857.9290) -- (238.0254,831.8217) -- (238.0254,831.8217) -- (235.4516,824.6953) -- (232.8931,818.8413) -- (231.4307,815.9062) -- (225.6891,806.4422) -- (223.8117,803.8898) -- (216.4890,795.6946) -- (211.6582,791.4400) -- (209.2911,789.5762) -- (209.2911,789.5762) -- (205.0664,784.4073) -- (200.6719,779.9106) -- (191.3599,772.4838) -- (190.5775,771.9609) -- (179.0758,765.5994) -- (170.8489,762.2332) -- (167.1783,760.9746) -- (141.0710,763.1719) -- (90.0567,768.3183) -- (9.5098,777.1664) -- (9.5098,989.4902);

        \path[fill=gray!80,draw=black,draw opacity=0.900,line cap=round,line join=round,line width=0.500pt,miter limit=10.00](989.6961,9.3039) -- (777.1664,9.5098) -- (777.1664,9.5098) -- (769.6814,77.6250) -- (764.2275,131.4150) -- (760.9746,167.1783) -- (760.9746,167.1783) -- (761.1063,167.5812) -- (761.1063,167.5812) -- (763.5484,174.3047) -- (766.1069,180.1587) -- (767.5693,183.0937) -- (773.3109,192.5578) -- (775.1883,195.1102) -- (782.5110,203.3054) -- (787.3418,207.5600) -- (789.4949,209.2305) -- (789.7089,209.4238) -- (789.7089,209.4238) -- (790.3502,210.2721) -- (790.3502,210.2721) -- (795.2782,216.0527) -- (798.3281,219.0894) -- (807.6401,226.5162) -- (808.4225,227.0391) -- (819.9242,233.4006) -- (828.1511,236.7668) -- (831.8217,238.0254) -- (835.5474,237.8911) -- (886.2188,233.0369) -- (956.5313,225.5233) -- (989.4902,221.8336) -- (989.6961,9.3039);

        \path[pattern= dots, draw=black,draw opacity=0.900,line cap=round,line join=round,line width=0.500pt,miter limit=10.00] (499.5000,989.6961) -- (499.5000,499.5000) -- (989.6961 ,499.5000) -- (989.6961 ,989.6961) -- (499.5000,989.6961);

        \path[pattern= dots, draw=black,draw opacity=0.900,line cap=round,line join=round,line width=0.500pt,miter limit=10.00]  (499.5000,499.5000) -- (499.5000, 9.3039) -- (9.3039 ,9.3039) -- (9.3039, 499.5000) -- (499.5000,499.5000);

        \path[pattern= dots, draw=black,draw opacity=0.900,line cap=round,line join=round,line width=0.500pt,miter limit=10.00] (126.9711,618.3268) -- (126.9711,618.3268) -- (51.2362,600.4115) -- (9.3039,590.7454) -- (9.3039,662.8987) -- (9.3039,662.8987) -- (41.2480,651.8920) -- (72.0098,640.5685) -- (101.5508,628.9560) -- (126.9711,618.3268);

        \path[pattern = dots, draw=black,draw opacity=0.900,line cap=round,line join=round,line width=0.500pt,miter limit=10.00] (380.6732,872.0289) -- (380.6732,872.0289) -- (396.2285,937.6325) -- (408.2546,989.6961) -- (336.1013,989.6961) -- (336.1013,989.6961) -- (345.6937,961.7396) -- (355.3931,935.0469) -- (365.2818,909.3307) -- (375.4795,884.2656) -- (380.6732,872.0289);

        \path[pattern = dots,draw=black,draw opacity=0.900,line cap=round,line join=round,line width=0.500pt,miter limit=10.00] (989.6961,408.2546) -- (989.6961,408.2546) -- (903.3901,388.1719) -- (872.0289,380.6732) -- (872.0289,380.6732) -- (872.0289,380.6732) -- (900.1348,368.9565) -- (929.4316,357.5051) -- (959.9492,346.3272) -- (989.6961,336.1013);

        \path[pattern = dots, draw=black,draw opacity=0.900,line cap=round,line join=round,line width=0.500pt,miter limit=10.00] (590.7454,9.3039) -- (590.7454,9.3039) -- (608.1426,84.2567) -- (618.3268,126.9711) -- (618.3268,126.9711) -- (618.3268,126.9711) -- (628.4927,102.6901) -- (638.5444,77.2995) -- (648.2889,51.2578) -- (657.8431,24.2396) -- (662.8987,9.3039);

        \path[pattern=dots, draw=black,draw opacity=0.900,line cap=round,line join=round,line width=0.500pt,miter limit=10.00] (224.2290,499.5000) -- (224.2290,499.5000) -- (231.6777,535.7173) -- (237.1550,563.9524) -- (237.1550,563.9524) -- (237.1550,563.9524) -- (260.7305,550.1834) -- (283.4355,536.0539) -- (305.6523,521.3202) -- (327.1367,506.1325) -- (336.1013,499.5000);

        \path[pattern=dots, draw=black,draw opacity=0.900,line cap=round,line join=round,line width=0.500pt,miter limit=10.00] (499.5000,774.7710) -- (499.5000,774.7710) -- (457.1118,766.1016) -- (435.0476,761.8450) --  (435.0476,761.8450) -- (435.0476,761.8450) -- (446.7654,741.6875) -- (458.9566,721.8307) -- (471.4125,702.6250) -- (484.3487,683.7448) -- (497.7874,665.1901) -- (499.5000,662.8987);

        \path[pattern = dots, draw=black,draw opacity=0.900,line cap=round,line join=round,line width=0.500pt,miter limit=10.00] (761.8450,435.0476) -- (761.8450,435.0476) -- (769.2754,473.0152) -- (774.7710,499.5000) -- (662.8987,499.5000) -- (662.8987,499.5000) -- (684.3145,483.9478) -- (706.2871,468.9844) -- (728.7480,454.6333) -- (751.9414,440.7252) -- (761.8450,435.0476);

        \path[pattern = dots, draw=black,draw opacity=0.900,line cap=round,line join=round,line width=0.500pt,miter limit=10.00] (563.9524,237.1550) -- (563.9524,237.1550) -- (520.7843,228.6667) -- (499.5000,224.2290) -- (499.5000,336.1013) -- (499.5000,336.1013) -- (512.8100,317.8594) -- (525.8453,298.9792) -- (538.3940,279.7734) -- (550.6736,259.9167) -- (562.6735,239.4089) -- (563.9524,237.1550);

\end{scope}

\end{tikzpicture}
   };

 \draw[->] (0,0) -- (0,2.3);
 \draw[->] (0,0) -- (2.3,0);
 
 \node at (5.1,-0.1) {1};
 \node at (-0.075,5.1) {1};
 \node at (-0.075,-0.075) {0};
 
\node[scale=1.2] at (2.3,-0.15) {$\xi$};
 \node[scale=1.2] at (-0.1,2.3) {$\eta$};
 
 \fill[fill=white] (0.2,1.89) circle (3pt);
    \node at (0.2,1.89) {$\mathcal{I}_2$};
    
 \fill[fill=white] (1.3,2.4) circle (3pt);
    \node at (1.3,2.4) {$\mathcal{I}_3$};    
    
\fill[fill=white] (1.25,3.75) circle (4.5pt);
    \node[scale=1.5] at (1.25,3.75) {$\mathcal{I}_1$};

\filldraw[pattern = dots] (5.5,3.3) rectangle (5.7,3.5);
\node at (6.43,3.38) {Elliptic islands};

\filldraw[fill=gray!80] (5.5,2.7) rectangle (5.7,2.9);
\node at (6.34,2.78) {Exp. mixing};

\node at (5.6,2.5) {$\subset$};

\filldraw[fill = gray!40] (5.5,2.1) rectangle (5.7,2.3);
\node at (6.12,2.18) {Mixing};

\node at (5.6,1.9) {$\subset$};

\filldraw[fill = gray!10] (5.5,1.5) rectangle (5.7,1.7);
\node at (6.28,1.58) {Hyperbolic};

    \end{tikzpicture}
    \caption{Proven behaviour of $H_{(\xi,\eta)}$ over the parameter space $0<\xi,\eta<1$. Regions $\mathcal{I}_n$ and their reflections exhibit persistent elliptic island structures, explored in section \ref{sec:islands}.}
    \label{fig:fullParameterSpace}
\end{figure}

The aim of this paper is to establish mixing results on some subset of the parameter space $0<\xi,\eta<1$. Figure \ref{fig:poincare} shows the orbit starting at $z_0=(1/\sqrt{2},1/\sqrt{3})$ for three parameter choices, highlighting a variety of ergodic and elliptic behaviour over the parameter space. A sketch of the results in the present work is given in Figure \ref{fig:fullParameterSpace}. In section \ref{sec:outline} we state two results which give an efficient scheme for proving the mixing property. Sections \ref{sec:lyp}, \ref{sec:man}, \ref{sec:repMan} give a simplified proof of the mixing property over a subset of the one dimensional parameter space $\xi = \eta$.  In particular we show:
\begin{thm}
\label{thm:Bernoulli}
The map $H = H_{(\eta,\eta)}$ has the Bernoulli property for $0<\eta<\eta_1 \approx 0.2024$.
\end{thm}
In section \ref{sec:rates} we prove exponential decay of correlations for $H$ over this subset. In section \ref{sec:param} we generalise our results to the wider two dimensional parameter space; making a small technical adjustment to Lemmas \ref{lemma:unstableGrowth} and \ref{lemma:stableGrowth}, establishing parameter space symmetries, and proving our main result, Theorem \ref{thm:2dMixingResults}. Section \ref{sec:otm} looks at the special case $H_{(\xi,\eta)}$ with $\xi=\eta=\frac{1}{2}$, which exhibits some unique dynamical features. We conclude with some final remarks in section \ref{sec:discussion}.

%% file: sections/outline.tex
Our scheme for proving the Bernoulli property is to satisfy the qualifications given in the following theorem from \cite{katok_invariant_1986}, paraphrased in \cite{sturman_ottino_wiggins_2006}.

\begin{thm}[\citeauthor{katok_invariant_1986}]
\label{thm:katok-strelcyn}
Let $f: X \rightarrow X$ be a measure preserving map on a measure space $(X,\mathcal{F},\mu)$ such that $f$ is $C^2$ smooth outside of a \emph{singularity set} $S$ where differentiability fails. Suppose that the Katok-Strelcyn conditions hold:
\begin{enumerate}[label={\bfseries (KS\arabic*):}]
    \item $\exists \, a,C_1>0$ s.t. $\forall \, \epsilon>0$, $\mu(B_\varepsilon(S)) \leq C_1 \varepsilon^a$.
    \item $\exists \, b,C_2>0$ s.t. $\forall \, z \in X \setminus S$, $||D^2_zf|| \leq C_2 \,d(z,S)^{-b}$ where $D^2_zf$ is the second derivative of $f$ at $z$.
    \item Lyapunov exponents exist and are non-zero almost everywhere.
\end{enumerate}
Then at almost every $z$ we can define local unstable and stable manifolds $\gamma_u(z)$ and $\gamma_s(z)$. Suppose that the manifold intersection property holds:
\begin{enumerate}[label={\bfseries (M):}]
    \item For almost any $z,z'\in X$, $\exists \, m,n$ s.t. $f^m(\gamma_u(z)) \cap f^{-n}(\gamma_s(z')) \neq \varnothing$.
\end{enumerate}
Then $f$ is ergodic. Furthermore the Bernoulli property holds, provided we can show the repeated manifold intersection property:
\begin{enumerate}[label={\bfseries (MR):}]
    \item For almost any $z,z'\in X$ there exists $M,N$ such that for all $m>M$ and $n>N$, $f^m(\gamma_u(z)) \cap f^{-n}(\gamma_s(z')) \neq \varnothing$.
\end{enumerate}
\end{thm}

The scheme extends Pesin theory (establishing ergodic properties of $C^2$ smooth non-uniformly hyperbolic systems, \citealp{pesin_characteristic_1977}) to systems which are smooth outside of some singularity set. The conditions \textbf{(KS1-2)} ensure that this set has manageable influence, and follow easily from our map's definition. We take our map as $f = H$, our domain as $X = \tor$, and our singularity set as $S = \mathcal{D}$. Taking $\mu$ to be the Lebesgue measure on $\tor$, clearly $\mu(S)=0$. When we say `for almost any $z \in \tor$', we will be referring to the full measure set $X' = \tor \setminus S_\infty$, $S_\infty = \bigcup_{k \geq 0} H^{-k}(\mathcal{D}) \cup \bigcup_{k \geq 0} H^{k}(\mathcal{D}')$, where $H$ and all its powers $H^k$, $k\in \mathbb{Z}$ are differentiable. Since we can cover $\mathcal{D}$ with arbitrarily thin rectangles, \textbf{(KS1)} follows for some $C_1>0$ with $a=1$. Since $H$ is piecewise linear, \textbf{(KS2)} follows trivially.

Moving onto \textbf{(KS3)}, we define the (forwards-time) Lyapunov exponent at a point $z \in \mathbb{T}^2$ in direction $v\in \mathbb{R}^2$ by
\[ \chi(z,v) = \limn \frac{1}{n} \log||DH^n_z v||, \]
where \[ DH^n_z = DH_{H^{n-1}(z)}\cdot ... \cdot DH_{H(z)} \cdot DH_z \]
is well defined at almost every $z$. We define $\log^+(\cdot) = \max \{ \log(\cdot), 0 \}$ and let $||\cdot||_{\mathrm{op}}$ be the operator norm. Existence of Lyapunov exponents almost everywhere follows from Oseledets' theorem (\citealp{oseledets_multiplicative_1968}) provided that $\log^+||DH||_{\mathrm{op}}$ is integrable. This clearly holds, so 
our first substantial task is proving that these exponents are non-zero. A particular form of Oseledets' theorem in two dimensions is useful here. We paraphrase from \cite{viana_lectures_2014}:

\begin{thm}[\citeauthor{oseledets_multiplicative_1968}, \citeauthor{viana_lectures_2014}]
\label{thm:Oseledets-2d}
Let $F:X \times \mathbb{R}^2 \rightarrow X \times \mathbb{R}^2$ be given by $F(x,v) = (f(x),A(x)v)$ for some measure preserving map $f$ on a 2-dimensional manifold $X$ and some measurable function $A: X \rightarrow \mathrm{GL}(2)$. Suppose $ \log^+|| A^{\pm 1} ||$ are integrable and define
\[ \lambda_+(x) = \limn \frac{1}{n} \log||A^n(x)||, \quad \lambda_-(x) = \limn \frac{1}{n} \log|| (A^n(x))^{-1}||^{-1}, \]
where $A^n(x) = A(f^{n-1}(x))\cdot ... \cdot A(f(x)) \cdot A(x)$. Then for almost every $x\in X$,
\begin{enumerate}
    \item either $\lambda_-(x) =\lambda_+(x)$ and
    \[ \limn \frac{1}{n} \log || A^n(x)v || = \lambda_\pm(x) \quad \forall v \in \mathbb{R}^2 \setminus \{ 0 \} \]
    \item or $\lambda_+(x) >\lambda_-(x)$ and there exists a vector line $E_x^s \subset \mathbb{R}^2$ such that
    \[ \limn \frac{1}{n} \log || A^n(x)v || = \begin{dcases}
    \lambda_-(x) & \text{for } v \in E_x^s \setminus \{0\}, \\
     \lambda_+(x) & \text{for } v \in \mathbb{R}^2 \setminus E_x^s.
     \end{dcases}
    \]
\end{enumerate}

\end{thm}

\begin{cor}
\label{cor:OurCorollary}
Further assuming that $A$ takes values in $\mathrm{SL}(2)$ gives $\lambda_-(x) = -\lambda_+(x)$. Hence if at some $x$ there exists $v_0 \in \mathbb{R}^2$ with $\lim_n \frac{1}{n} \log || A^n(x)v_0 || \neq 0$, it follows that $\lim_n \frac{1}{n} \log || A^n(x)v || \neq 0$ for all non-zero vectors $v$.
\end{cor}

Applying this corollary to the cocycle generated by the derivative $DH$ of our map $H$ gives an efficient scheme for establishing non-zero Lyapunov exponents. We let $A^n(z) = DH_z^n$, which takes values in SL(2). If there exists $v_0$ such that $|| DH^n_z v_0 ||$ grows exponentially with $n$, Corollary \ref{cor:OurCorollary} gives $\chi(z,v) \neq 0$ for all $v \neq 0$.

%% file: sections/lyp.tex
\begin{prop}
\label{prop:nonZeroChi}
Let $0<\eta<1/2$. At almost every $z$, $\chi(z,v) \neq 0$ for all $v\neq 0$.
\end{prop}

\begin{proof}
Let $H = H_{(\eta,\eta)}$ and take matrices $M_j = DH|_{A_j}$, all of which are hyperbolic over $0<\eta<1/2$. Write the gradients of the unstable, stable eigenvectors of $M_j$ as $g_j^u$, $g_j^s$. One can verify that
\begin{equation}
    \label{eq:coneInequality}
     g_4^u(\eta) < g_2^u(\eta) < g_2^s(\eta) < g_1^s(\eta) < g_4^s(\eta) < g_3^s(\eta) < g_3^u(\eta) < g_1^u(\eta) 
\end{equation}
across $0<\eta<1/2$. This allows us to define a cone region $\mathcal{C}$ in the tangent space, bounded by the unstable eigenvectors of $M_2$ and $M_3$, which includes all of the unstable eigenvectors of the $M_j$, and none of the stable eigenvectors. It follows that this cone is \emph{invariant}, that is, $M_j \mathcal{C} \subset \mathcal{C}$ for each $j$. It is easily verified that this cone is \emph{expanding}, that is, there exists $\delta>0$ such that $|| M_j v  || \geq (1+\delta) ||v||$ for each $j$, vector $v\in \mathcal{C}$, where $|| \cdot ||$ is whatever norm we put on the tangent space. Lower bounds on these expansion factors using a convenient norm, $||\cdot||_\infty$, are given explicitly in the next section.

For any $z\in X'$, $v_0 \in \mathcal{C}$, it follows that $|| DH_z^n v_0 ||$ grows exponentially with $n$. By Corollary \ref{cor:OurCorollary}, this implies $\chi(z,v)\neq 0$ for any $v\neq 0$.
\end{proof}

Hence we have non-zero Lyapunov exponents almost everywhere, i.e. $H$ is hyperbolic.

%% file: sections/man.tex
In this section we will prove the following:
\begin{prop}
\label{prop:M}
Condition \textbf{(M)} holds for $H$ when $0<\eta<\eta_1\approx 0.2024$.
\end{prop}

The proof consists of three stages. Non-zero Lyapunov exponents at $z \in X'$ implies the existence of local unstable and stable manifolds $\gamma_u(z)$ and $\gamma_s(z)$ at $z$. The first stage, Lemma \ref{lemma:alignment}, describes the nature of these local manifolds. In the next stage, Lemmas \ref{lemma:unstableGrowth} and \ref{lemma:stableGrowth}, we give an iterative scheme for growing the backwards (forwards) images of any local (un)stable manifold. We then grow the images of these manifolds up until the point where the images connect up certain partition boundaries. This then allows us, by Lemmas \ref{lemma:vSegs} and \ref{lemma:hSegs}, to establish an intersection in the next several iterates.

Let $\mathcal{C}'$ be the cone bounded by the stable eigenvectors of $M_2$ and $M_3$, including the stable eigenvectors of each of the $M_j$. It follows that this cone is invariant and expanding under $H^{-1}$. The cones $\mathcal{C}$ and $\mathcal{C}'$ provide bounds on the gradients of local manifolds:

\begin{lemma}
\label{lemma:alignment}
At every $z \in X'$, $\gamma_u(z)$ is a line segment aligned with some $v\in \mathcal{C}$, and $\gamma_s(z)$ is a line segment aligned with some $v' \in \mathcal{C}'$.
\end{lemma}

\begin{proof}
Since $H$ is piecewise linear, $\gamma_u(z)$ and $\gamma_s(z)$ are line segments, aligned with some vectors $v$ and $v'$ respectively. By definition, for any $\zeta,\zeta' \in \gamma_u(z)$
\begin{equation}
    \label{eq:unstableManifold}
    \mathrm{dist}(H^{-n}(\zeta),H^{-n}(\zeta')) \rightarrow 0
\end{equation} 
as $n \rightarrow \infty$. Similarly for any $\zeta,\zeta' \in \gamma_s(z)$
\begin{equation}
    \label{eq:stableManifold}
    \mathrm{dist}(H^{n}(\zeta),H^{n}(\zeta')) \rightarrow 0
\end{equation}
as $n\rightarrow \infty$. Note that $\mathcal{C}$ can be described at the cone region bounded by the stable eigenvectors of $M_2^{-1}$ and $M_3^{-1}$ and including the stable eigenvectors of each $M_j^{-1}$. Clearly $v$ must be aligned in this cone, for if it falls outside of this region, repeatedly applying the $M_j^{-1}$ will pull $v$ towards the invariant expanding cone $\mathcal{C}'$, resulting in exponential growth in its norm, which contradicts (\ref{eq:unstableManifold}). Similarly $v'$ must lie in $\mathcal{C}'$ to avoid contradicting (\ref{eq:stableManifold}).

\end{proof}

We now move onto the growth stage. We say that a line segment has \emph{simple intersection} with $A_j$ if its restriction to $A_j$ is empty or a single line segment. We define the $x$- and $y$-diameters of a line segment $\Gamma$ by $\mathrm{diam}_x(\Gamma) = \nu \left( \{x \,|\, (x,y) \in \Gamma \} \right)$ and $\mathrm{diam}_y(\Gamma) = \nu \left( \{y \,|\, (x,y) \in \Gamma \} \right)$, where $\nu$ is the Lebesgue measure on $\mathbb{R}$.

\begin{lemma}[Growth Lemma]
\label{lemma:unstableGrowth}
Let $\eta < \eta_1$. Given a line segment $\Gamma_{p-1}$, aligned with some $v \in \mathcal{C}$ and having simple intersection with each $A_j$, there exists a line segment $\Gamma_p \subset H(\Gamma_{p-1})$ such that
\begin{enumerate}[label={(C\arabic*)}]
    \item $\Gamma_p$ is aligned with some vector in $\mathcal{C}$,
    \item $\mathrm{diam}_y(\Gamma_p) \geq (1+\delta)\,\mathrm{diam}_y(\Gamma_{p-1})$ for some $\delta >0$.
\end{enumerate}
\end{lemma}

\begin{proof}
Let $||\cdot||$ denote the $||\cdot||_\infty$ norm. Since $|v_2| \geq |v_1|$ for every $v = (v_1,v_2)^T \in \mathcal{C}$, vectors $v \in \mathcal{C}$ have norm $||v|| = |v_2|$. Define minimum diameter expansion factors
\[ K_j(\eta) = \inf_{v \in \mathcal{C}} \frac{|| M_j v ||}{|| v ||} \]
for each of the matrices $M_j$. Over the cone, $M_1$, $M_2$ are minimal on the cone boundary given by the unstable eigenvector $v_u(M_2)$ of $M_2$, and $M_3$, $M_4$ are minimal on the other cone boundary $v_u(M_3)$. From this we can calculate the $K_j$ as
\[ K_1 = \frac{2-\eta}{1-\eta},\quad K_2 = \frac{1-\eta}{\eta},\quad K_3 = \frac{1-\eta}{\eta},\text{ and } K_4 = \frac{1-\eta+\eta^2}{\eta^2}.\]

Suppose $\Gamma_{p-1}$ has non-simple intersection with all the $A_j$ and each intersection is non-empty. Write the restriction of $\Gamma_{p-1}$ to $A_j$ as $\Gamma^j$. Now if for some $j$
\begin{equation}
\label{eq:Kj}
K_j(\eta)\, \mathrm{diam}_y(\Gamma^j) >\mathrm{diam}_y(\Gamma_{p-1}),
\end{equation} 
we can take $\Gamma_p = H(\Gamma^j)$ to satisfy (C2). If $\Gamma_{p-1}$ was aligned with $v \in \mathcal{C}$, $\Gamma_{p-1}$ is now aligned with $M_j v \in \mathcal{C}$, so (C1) is also satisfied. If (\ref{eq:Kj}) does not hold, the proportion of $\Gamma^j$ in $\Gamma_{p-1}$ is bounded above by $K_j^{-1}$. Suppose (\ref{eq:Kj}) does not hold for $j=2,3,4$. Then the proportion of $\Gamma^1$ in $\Gamma_{p-1}$ is bounded below by
\[ \frac{\mathrm{diam}_y(\Gamma^1)}{\mathrm{diam}_y(\Gamma_{p-1})} > 1-\frac{1}{K_2} -\frac{1}{K_3} - \frac{1}{K_4}. \]
Hence taking $\Gamma_p = H(\Gamma^1)$ satisfies (C2) provided that
\begin{equation}
    \label{eq:expansionInequality}
    K_1(\eta) > \frac{1}{1-\frac{1}{K_2} -\frac{1}{K_3} - \frac{1}{K_4}}.
\end{equation}
Plugging in the expressions for $K_j(\eta)$, the above is satisfied for $0<\eta<\eta_1$, where $\eta_1 \approx 0.2024$ is the smallest real solution to the quartic equation $-\eta^4 +4\eta^3 -6\eta^2 + 6\eta -1  = 0$. So for $\eta$ in this range, choosing one of $\Gamma_p = H(\Gamma^j)$ will always satisfy (C2). The case where $\Gamma_{p-1}$ has empty intersection with one or more of the $A_j$ follows as a trivial consequence.
\end{proof}

The equivalent lemma for the growth of line segments under $H^{-1}$ is as follows. Recall the partition of the torus into four sets $A_j'$ given in Figure \ref{fig:Apartition}.

\begin{lemma}
\label{lemma:stableGrowth}
Let $\eta < \eta_1$. Given a line segment $\Gamma_{p-1}$, aligned with some $v' \in \mathcal{C}'$ and having simple intersection with each $A_j'$, there exists a line segment $\Gamma_p \subset H^{-1}(\Gamma_{p-1})$ such that
\begin{enumerate}[label={(C\arabic*')}]
    \item $\Gamma_p$ aligned with some vector in $\mathcal{C}'$,
    \item $\mathrm{diam}_x(\Gamma_p) \geq (1+\delta)\,\mathrm{diam}_x(\Gamma_{p-1})$ for some $\delta >0$.
\end{enumerate}
\end{lemma}

\begin{proof}
Argument is entirely analogous. It is easily verified that the minimum diameter expansion of $M_j^{-1}$ over $\mathcal{C}'$ is $K_j$ as defined before. The lemma, then, also holds for $0<\eta<\eta_1$.
\end{proof}

Moving onto the final mapping stage, call any line segment $\Gamma \subset R_1$ which joins the upper and lower boundaries ($y=0$, $y=1-\eta$) a \emph{$v$-segment}. Similarly we call any line segment $\Gamma \subset R_1$ which joins the left and right boundaries ($x=0$, $x=1-\eta$) a \emph{$h$-segment}. Clearly $v$-segments and $h$-segments must always intersect.

\begin{lemma}[Mapping Lemma]
\label{lemma:vSegs}
Let $\Gamma$ be a line segment contained within some $A_j'$. If $\Gamma$ has non-simple intersection with some $A_j$, then $H^k(\Gamma)$ contains a $v$-segment for some $k \in \{ 1,2,3,4 \}$.
\end{lemma}

\begin{proof}

Note that the sets $A_1'$, $A_3'$ are entirely contained within the strip $\{ x \leq 1-\eta \}$, and the sets $A_2'$, $A_4'$ are entirely contained within the strip $\{ x \geq 1-\eta \}$, so $\Gamma$ lies entirely within one of these strips. Suppose first that it lies in $\{ x \leq 1-\eta \}$, then $\Gamma$ must have non-simple intersection with $A_1$ or $A_3$. Non-simple intersection with $A_2$ and $A_4$ is possible, but involves wrapping vertically around the torus, and in doing so implies non-simple intersection with $A_1$ or $A_3$. Assume $\Gamma$ has non-simple intersection with $A_1$. Then it must either connect the segments 2a and 2b (shown in Figure \ref{fig:lowerCaseSegs}) though $A_2$ or connect the segments 4a and 4b through $A_4$, depending which way it connects the two parts of $A_1$. The same is true when $\Gamma$ has non-simple intersection with $A_3$.

Equivalent analysis can be applied to the strip $\{x\geq 1-\eta \}$. For $\Gamma$ in this strip, it follows that $\Gamma$ connects 3a to 3b through $A_3$ or connects 1a to 1b through $A_1$. This gives four possible cases. Denote the case where $\Gamma$ connects $j$a to $j$b through $A_j$ by case ($j$). We will show that all cases reduce to case (3). Suppose first that $\Gamma$ satisfies case (4), connecting 4a to 4b through $A_4$. Then $H(\Gamma)$ connects 4a' to 4b' through $A_4'$ (see Figure \ref{fig:lowerCaseSegs}). To do this, $H(\Gamma)$ must connect the segments 1a and 1b, passing through $A_1$. That is, $H(\Gamma)$ satisfies case (1). One can similarly show that if $\Gamma$ satisfies case (1) then $H(\Gamma)$ satisfies case (2), and that if $\Gamma$ satisfies case (2) then $H(\Gamma)$ satisfies case (3).

Looking at the images $\mathrm{3a'} = H(\mathrm{3a})$ and $\mathrm{3b'} = H(\mathrm{3b})$, we see that any line segment joining 3a' to 3b' must pass through $y=0$ and $y=1-\eta$, the lower and upper boundaries of $R_1$. It follows that if $\Gamma$ satisfies case (3), $H(\Gamma)$ contains a $v$-segment. For any of the four cases ($j$), $j=1,2,3,4$, $H^k(\Gamma)$ will contain a $v$-segment for $k=3,2,1,4$.
\begin{figure}
    \centering

     \begin{tikzpicture}
     \tikzmath{\e = 0.25;} 
    
    \node[scale=1.7] at (-4,0) {
    \begin{tikzpicture}
    
    \fill[fill=gray!40] (0, {4*(1-\e+\e*\e)}) -- ({4*(1-\e)},4) -- (0,4) --  (0, {4*(1-\e+\e*\e)});
    
    \fill[fill=gray!40] (4,4) -- (0,{4*(1-\e)}) -- ({4*(1-\e)},{4*(1-\e)}) -- (4, {4*(1-\e+\e*\e)}) -- (4,4);
    
    \fill[fill=gray!20] ({4*(1-\e)},{4*(1-\e)}) -- (4,{4*(1-\e)}) -- (4, {4*(1-\e)^2}) -- ({4*(1-\e)},{4*(1-\e)});
    
    \fill[fill=gray!20] (0,{4*(1-\e)}) -- (4,0) -- ({4*(1-\e)},0) -- (0, {4*(1-\e)^2}) -- (0,{4*(1-\e)});

    \node at (2-0.5*4*\e,2-0.5*4*\e) {$A_2$};
    
    \node at (4-0.5*4*\e,2-0.5*4*\e) {$A_1$};

    \node at (4-0.5*4*\e,4-0.5*4*\e) {$A_3$};
    
    
    \node[scale=0.6] (a) at ({0.5*4*(1-\e)},{0.5*(4-4*\e*\e)+0.5*4*(1-\e)}) {4a};
    \draw [thick, -|] (0,{4*(1-\e)})  -- (a) --  ({4*(1-\e)},4-4*\e*\e);
    
    \node[scale=0.6] (a) at ({0.5*4*(1-\e)},{0.5*4*(1-\e+\e*\e)+0.5*4}) {4b};
    \draw [thick] (0, {4*(1-\e+\e*\e)})  -- (a) --  ({4*(1-\e)},4);

    
    \node[scale=0.6] (a) at ({0.5*4    +0.5*4*(1-\e)  },{0.5*4*(1-\e)^2   +0.5*4*(1-\e) }) {1b};
    \draw [thick] (4, {4*(1-\e)^2})  -- (a) --  ({4*(1-\e)},{4*(1-\e)});

    \node[scale=0.6] (a) at ({0.5*4*(1-\e)    +0.5*4    },{0.5* 4*(1-\e)*\e  +0.5*0  }) {1a};
    \draw [thick,-|]  (4,0) -- (a) -- ({4*(1-\e)}, {4*(1-\e)*\e}) ;
    
    \node[scale=0.6] (a) at ({0.5*4*(1-\e)    +0.5*0    },{0.5*4*(1-\e)*\e   +0.5*4*(1-\e)  }) {2a};
    \draw [thick] ({4*(1-\e)}, {4*(1-\e)*\e})  -- (a) --  (0,{4*(1-\e)});
 
    \node[scale=0.6] (a) at ({0.5*0    +0.5*4*(1-\e)    },{0.5*4*(1-\e)^2   +0.5*0  }) {2b};
    \draw [thick] (0, {4*(1-\e)^2})  -- (a) --  ({4*(1-\e)},0);
    
    \node[scale=0.6] (a) at ({0.5*4*(1-\e)    +0.5*4    },{0.5*(4-4*\e*\e)   +0.5*4  }) {3a};
    \draw [thick] ({4*(1-\e)},4-4*\e*\e)  -- (a) --  (4,4);
    
    \node[scale=0.6] (a) at ({0.5*4*(1-\e)    +0.5*4    },{0.5*4*(1-\e)   +0.5*4*(1-\e+\e*\e) }) {3b};
    \draw [thick] ( {4*(1-\e)}, {4*(1-\e)})  -- (a) --  (4 ,{4*(1-\e+\e*\e)}  );

    \end{tikzpicture}
    };
    
    \node[scale=1.7] at (4.5,0) {
    \begin{tikzpicture}
    
    \fill[fill = gray!40] (0,4-4*\e) -- ({4*\e*(1-\e)},4) -- (0,4) -- (0,4-4*\e);
    \fill[fill = gray!40] (0,0) -- (4-4*\e,4) -- (4-4*\e,4-4*\e) -- ({4*\e*(1-\e)},0) -- (0,0);
    
    \fill[fill=gray!20] (4-4*\e,0) -- (4-4*\e,4-4*\e) -- (4-4*\e*\e,0) -- (4-4*\e,0);
    \fill[fill=gray!20] (4,0) -- (4-4*\e,4) -- (4-4*\e*\e,4)-- (4,4-4*\e) -- (4,0);

    
    \node[scale=0.6] (a) at (0,{0.5*4*\e*(1-\e)   +0.5*0  }) {1a'};
    \draw [thick] ( 0, 0)  -- (a) --  (0 , {4*\e*(1-\e)} );

    \node[scale=0.6] (a) at (-0.2,{0.5*4*(1-\e*\e)   +0.5*4  }) {3a'};
    \draw [thick] (0, 4) -- (0 , {4*(1-\e*\e)} );
    
    \node[scale=0.6] (a) at (4,{0.5*4*\e*(1-\e)   +0.5*4*(1-\e)  }) {2a'};
    \draw [thick] (4, {4*\e*(1-\e)}  )  -- (a) --  (4 , {4*(1-\e)});
    
    \node[scale=0.6] (a) at (4,{0.5*4*(1-\e*\e)   +0.5*4*(1-\e)  }) {4a'};
    \draw [thick, -|] (4, {4*(1-\e*\e)}  )  -- (a) --  (4 , {4*(1-\e)});
    
    \node[scale=0.6] (a) at (4-4*\e-0.2,{0.5*4*(1-\e+\e*\e)   +0.5*4*(1-\e)  }) {3b'};
    \draw [thick, -|] (4-4*\e, {4*(1-\e+\e*\e)}  )  --  (4-4*\e , {4*(1-\e)});

    \node[scale=0.6] (a) at (4-4*\e,{0.5*4*(1-\e+\e*\e)   +0.5*4  }) {4b'};
    \draw [thick, -|]  (4-4*\e , {4})   -- (a) --  (4-4*\e, {4*(1-\e+\e*\e)}  );
    
    \node[scale=0.6] (a) at (4-4*\e,{0.5*4*(1-\e)*(1-\e)   +0.5*4*(1-\e)  }) {1b'};
    \draw [thick, -|]  (4-4*\e , {4*(1-\e)})  -- (a) -- (4-4*\e, {4*(1-\e)*(1-\e)} );
    
    \node[scale=0.6] (a) at (4-4*\e,{0.5*4*(1-\e)*(1-\e)   +0.5*0  }) {2b'};
    \draw [thick]  (4-4*\e , {0})  -- (a) -- (4-4*\e, {4*(1-\e)*(1-\e)} );

    \end{tikzpicture}
    };
    \draw[->] (0,0) -- (1,0);
    \node[scale=1.5] at (0.5,0.5) {$H$};

    \end{tikzpicture}
    \caption{Four pairs of line segments $j\mathrm{a}$, $j\mathrm{b}$ on the boundaries of $A_j$, and their images $j\mathrm{a}'$, $j\mathrm{b}'$ under $H$ on the boundaries of $A_j'$.}
    \label{fig:lowerCaseSegs}
\end{figure}
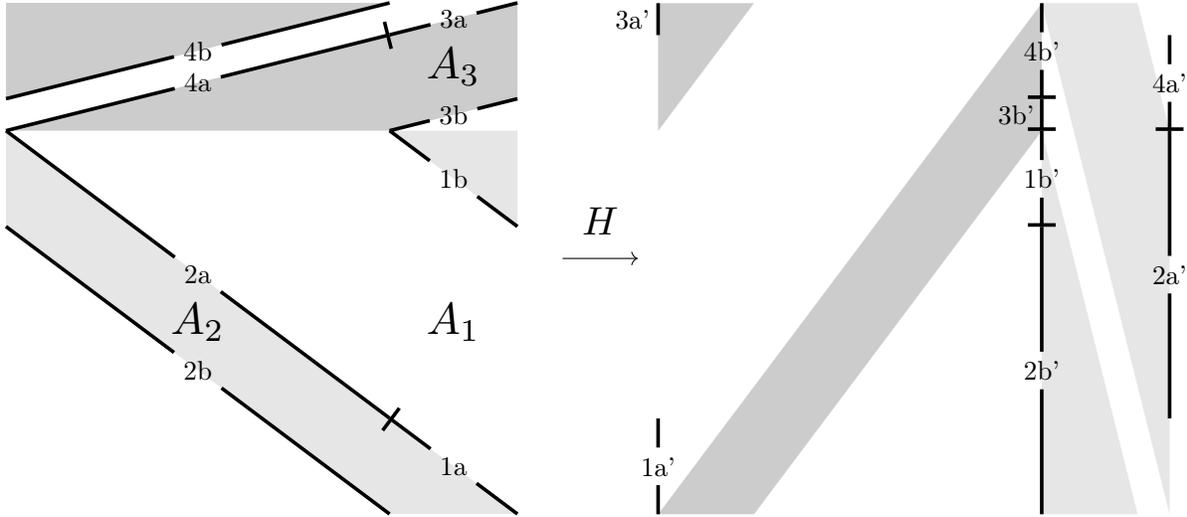
\end{proof}

\begin{lemma}
\label{lemma:hSegs}
Let $\Gamma$ be a line segment contained within some $A_j$. If $\Gamma$ has non-simple intersection with some $A_j'$, then $H^{-k}(\Gamma)$ contains a $h$-segment for some $k \in \{ 1,2,3,4 \}$.
\end{lemma}

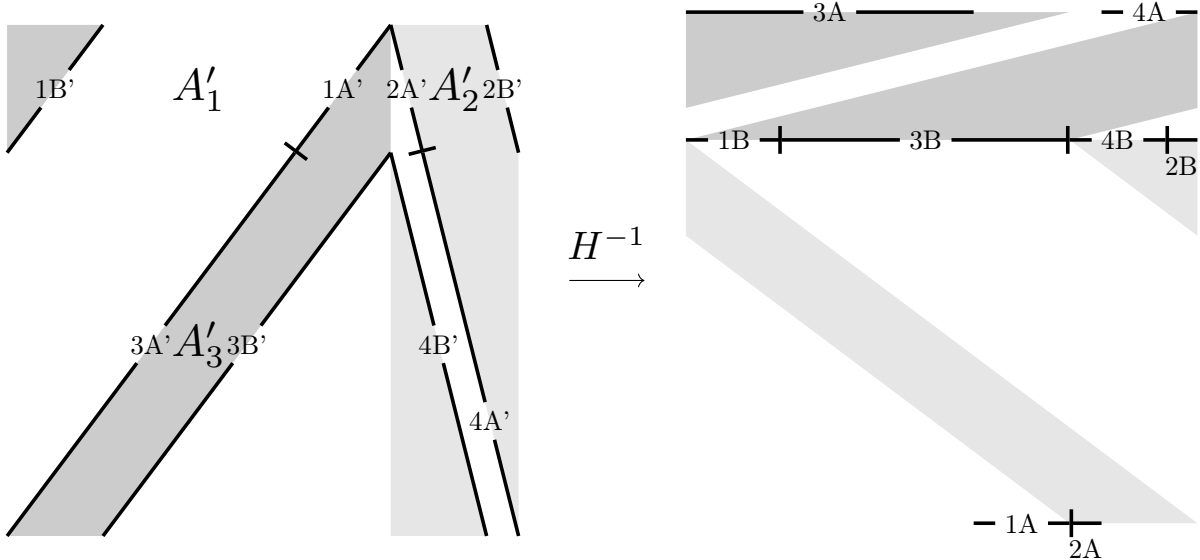
\begin{figure}
    \centering
    \begin{tikzpicture}

     \tikzmath{\e = 0.25;} 
     \node[scale=1.7] at (-4,0) {
    \begin{tikzpicture}
    \fill[fill = gray!40] (0,4-4*\e) -- ({4*\e*(1-\e)},4) -- (0,4) -- (0,4-4*\e);
    \fill[fill = gray!40] (0,0) -- (4-4*\e,4) -- (4-4*\e,4-4*\e) -- ({4*\e*(1-\e)},0) -- (0,0);
    
    \fill[fill=gray!20] (4-4*\e,0) -- (4-4*\e,4-4*\e) -- (4-4*\e*\e,0) -- (4-4*\e,0);
    \fill[fill=gray!20] (4,0) -- (4-4*\e,4) -- (4-4*\e*\e,4)-- (4,4-4*\e) -- (4,0);
    
    \node[scale=0.6] (a) at ({0.5*0    +0.5*4*\e*(1-\e) },{0.5*4*(1-\e) +0.5*4  }) {1B'};
    \draw [thick] (0 , {4*(1-\e)})  -- (a) --  ({4*\e*(1-\e)} , 4 );
    
    \node[scale=0.6] (a) at ({0.5*4*(1-\e)*(1-\e)    +0.5*4*(1-\e) },{0.5*4*(1-\e) +0.5*4  }) {1A'};
    \draw [thick, -|] ({4*(1-\e)} , 4 )  -- (a) --    ( {4*(1-\e)*(1-\e) }, {4*(1-\e)}) ;
    
    \node[scale=0.6] (a) at ({0.5*4*(1-\e)*(1-\e)    +0.5*4*0 },{0.5*4*(1-\e) +0.5*0  }) {3A'};
    \draw [thick] (0 , 0 )  -- (a) --    ( {4*(1-\e)*(1-\e) }, {4*(1-\e)}) ;
    
    \node[scale=0.6] (a) at ({0.5*4*(1-\e)   +0.5*4*\e*(1-\e) },{0.5*4*(1-\e) +0.5*0  }) {3B'};
    \draw [thick] ( {4*(1-\e)} , {4*(1-\e)} )  -- (a) --    ( {4*(1-\e)*\e }, 0) ;
    
    \node[scale=0.6] (a) at ({0.5*4*(1-\e)   +0.5*4*(1-\e*\e) },{0.5*4*(1-\e) +0.5*0  }) {4B'};
    \draw [thick] ( {4*(1-\e)} , {4*(1-\e)} )  -- (a) --    ( {4*(1-\e*\e) }, 0) ;
    
    \node[scale=0.6] (a) at ({0.3*4*(1-\e+\e*\e)    +0.7*4 },{0.3*4*(1-\e) +0.7*0  }) {4A'};
    \draw [thick] (4 , 0 )  -- (a) --    ( {4*(1-\e+\e*\e) }, {4*(1-\e)}) ;
    
    \node[scale=0.6] (a) at ({0.5*4*(1-\e+\e*\e)    +0.5*4*(1-\e) },{0.5*4*(1-\e) +0.5*4  }) {2A'};
    \draw [thick, -|] ({4*(1-\e)} , 4 )  -- (a) --    ( {4*(1-\e+\e*\e) }, {4*(1-\e)}) ;
    
    \node[scale=0.6] (a) at ({0.5*4*(1-\e*\e)    +0.5*4 },{0.5*4 +0.5*4*(1-\e)  }) {2B'};
    \draw [thick] ({4*(1-\e*\e)} , 4 )  -- (a) --    ( {4}, {4*(1-\e)}) ;

    \node at (2-0.5*4*\e,2-0.5*4*\e) {$A_3'$};
    
    \node at (1.5*4*\e,4-0.5*4*\e) {$A_1'$};
    \node at (4-0.5*4*\e,4-0.5*4*\e) {$A_2'$};

    \end{tikzpicture}};
    
    \draw[->] (0,0) -- (1,0);
    \node[scale=1.5] at (0.5,0.5) {$H^{-1}$};
    
    \node[scale=1.7] at (5,0) {
    \begin{tikzpicture}
    
     \fill[fill=gray!40] (0, {4*(1-\e+\e*\e)}) -- ({4*(1-\e)},4) -- (0,4) --  (0, {4*(1-\e+\e*\e)});
    
    \fill[fill=gray!40] (4,4) -- (0,{4*(1-\e)}) -- ({4*(1-\e)},{4*(1-\e)}) -- (4, {4*(1-\e+\e*\e)}) -- (4,4);
    
    \fill[fill=gray!20] ({4*(1-\e)},{4*(1-\e)}) -- (4,{4*(1-\e)}) -- (4, {4*(1-\e)^2}) -- ({4*(1-\e)},{4*(1-\e)});
    
    \fill[fill=gray!20] (0,{4*(1-\e)}) -- (4,0) -- ({4*(1-\e)},0) -- (0, {4*(1-\e)^2}) -- (0,{4*(1-\e)});

    \node[scale=0.6] (a) at ({0.5*0    +0.5*4*\e*(1-\e) },{0.5*4*(1-\e) +0.5*4*(1-\e)  }) {1B};
    \draw [thick, -|] (0 , {4*(1-\e)})  -- (a) --  ({4*\e*(1-\e)} , {4*(1-\e)} );
    
    \node[scale=0.6] (a) at ({0.5*4*(1-\e)*(1-\e)    +0.5*4*(1-\e) },{0.5*0 +0.5*0  }) {1A};
    \draw [thick] ({4*(1-\e)} , 0 )  -- (a) --    ( {4*(1-\e)*(1-\e) }, {0}) ;
    
    \node[scale=0.6] (a) at ({0.5*4*(1-\e)*(1-\e)    +0.5*4*0 },{0.5*4 +0.5*4  }) {3A};
    \draw [thick] (0 , 4 )  -- (a) --    ( {4*(1-\e)*(1-\e) }, 4) ;
    
    \node[scale=0.6] (a) at ({0.5*4*(1-\e)   +0.5*4*\e*(1-\e) },{0.5*4*(1-\e) +0.5*4*(1-\e)  }) {3B};
    \draw [thick, -|]   ( {4*(1-\e)*\e },{4*(1-\e)} ) -- (a) --    ( {4*(1-\e)} , {4*(1-\e)} );
    
    \node[scale=0.6] (a) at ({0.5*4*(1-\e)   +0.5*4*(1-\e*\e) },{0.5*4*(1-\e) +0.5*4*(1-\e)  }) {4B};
    \draw [thick] ( {4*(1-\e)} , {4*(1-\e)} )  -- (a) --    ( {4*(1-\e*\e) }, {4*(1-\e)}) ;
    
    \node[scale=0.6] (a) at ({0.5*4*(1-\e+\e*\e)    +0.5*4 },{4}) {4A};
    \draw [thick] (4 , 4 )  -- (a) --    ( {4*(1-\e+\e*\e) }, 4) ;
    
    \node[scale=0.6] (a) at ({0.5*4*(1-\e+\e*\e)    +0.5*4*(1-\e) },{-0.2}) {2A};
    \draw [thick, -|]    ( {4*(1-\e+\e*\e) }, 0) -- ({4*(1-\e)} , 0) ;
    
    \node[scale=0.6] (a) at ({0.5*4*(1-\e*\e)    +0.5*4 },{0.5*4*(1-\e) +0.5*4*(1-\e)  -0.2 }) {2B};
    \draw [thick, -|]  ( {4}, {4*(1-\e)}) -- ({4*(1-\e*\e)} , {4*(1-\e)} );
    
    \end{tikzpicture}};
    
    \end{tikzpicture}

    \caption{Four pairs of line segments $j\mathrm{A}'$, $j\mathrm{B}'$ on the boundaries of $A_j'$, and their images $j\mathrm{A}$, $j\mathrm{B}$ under $H^{-1}$ on the boundaries of $A_j$.}
    \label{fig:upperCaseSegs}
\end{figure}

\begin{proof}
The argument is almost entirely analogous, we say that $\Gamma$ satisfies case ($j'$) if $\Gamma$ connects $j$A' to $j$B' through $A_j'$ (see Figure \ref{fig:upperCaseSegs} for an illustration of the relevant segments). $\Gamma$ is entirely contained within one of the strips $\{y\leq 1-\eta\}$ or $\{y\geq 1-\eta\}$ which, together with the fact that $\Gamma$ has non-simple intersection with some $A_j'$, implies that $\Gamma$ satisfies case ($j$') for some $j$.
Again, we have that if $\Gamma$ satisfies case (4') then $H^{-1}(\Gamma)$ satisfies case (1'). This reduces to case (3'), and in turn reduces to case (2'). Any segment connecting 2A to 2B through $A_2$ must pass through the lines $x=1-\eta$ and $x=0$, the right and left boundaries of $R_1$. It follows that for $\Gamma$ satisfying case ($j$'), $j=1,2,3,4$, $H^{-k}(\Gamma)$ contains a $h$-segment for $k= 3,1,2,4$.
\end{proof}

We are now ready establish ergodicity.

\begin{proof}[Proof of Proposition \ref{prop:M}]
Given $z \in X'$, by Lemma \ref{lemma:alignment}, $\Gamma_0 = \gamma_u(z)$ is a line segment aligned with some vector $v \in \mathcal{C}$. By Lemma \ref{lemma:unstableGrowth} we can generate a sequence of line segments $(\Gamma_p)_{0\leq p \leq P}$, with $\Gamma_p \subset H^p(\gamma_u(z))$ and the diameter of $\Gamma_p$ growing exponentially with $p$. It follows that after finitely many $P$ steps, $\Gamma_P$ must have non-simple intersection with one of the partition elements $A_j$. Since $H^{-1}(\Gamma_P)$ lies entirely within some $A_j$, $\Gamma_P$ lies entirely within some $A_j'$. Now by Lemma \ref{lemma:vSegs}, $H^k(\Gamma_P)$ contains a $v$-segment for some $k \in \{1,2,3,4\}$. Hence we have found $m = P + k$ such that $H^m(\gamma_u(z))$ contains a $v$-segment. Similarly given $z'\in X'$, by Lemmas \ref{lemma:alignment}, \ref{lemma:stableGrowth}, and \ref{lemma:hSegs}, we can find $n$ such that $H^{-n}(\gamma_s(z'))$ contains a $h$-segment. It follows that they must intersect.
\end{proof}

We now move onto establishing stronger mixing properties.

%% file: sections/repMan.tex
\begin{prop}
\label{prop:MR}
Condition \textbf{(MR)} holds for $H$ when $0<\eta<\eta_1\approx 0.2024$.
\end{prop}

\begin{proof}

Given $z\in X'$, by Lemmas \ref{lemma:alignment}, \ref{lemma:unstableGrowth}, \ref{lemma:vSegs}, we have found $M_0$ such that $H^{M_0}(\gamma_u(z))$ contains a segment $\Gamma$ which joins 3a' to 3b' through $A_3'$. As shown in the previous section, this means that $\Gamma$ contains a $v$-segment. It also follows that $\Gamma$ satisfies case (2) so, by induction, we have that $H^{2k}(\Gamma)$ contains a $v$-segment for $k \in \mathbb{N}$. Consider the quadrilateral $Q_1\subset A_1$, defined by its corners
\[ q_1 = \left( \frac{(1-\eta)^3}{1+(1-\eta)^2},0 \right), \quad q_2 = \left(\frac{(1-\eta)^2}{1+(1-\eta)^2},0 \right), \quad q_3 = \left(0,\frac{(1-\eta)^3}{1+2(1-\eta)^2} \right), \quad q_1 = \left(0,\frac{(1-\eta)^4}{1+2(1-\eta)^2} \right). \]
An illustration of $Q_1$ and its image $H(Q_1) \subset A_2$ are shown in Figure \ref{fig:rmip}. One can show that each of the points $q_i$ map into the boundary of $A_2$ so that if $\Gamma$ joins the dashed boundaries of $Q_1$, then $H(\Gamma)$ satisfies case (2). For $\Gamma$ joining 3a' to 3b' through $A_3'$, $\Gamma$ must intersect the line $y=0$ at some point $(x,0)$ with $0 \leq x \leq x_v = \eta^2(1-\eta)\left(1-\eta+\eta^2\right)^{-1}$. Hence our $\Gamma$ joins the dashed lines of $Q_1$ as described, provided that $x_v(\eta)\leq q_1(\eta)$. This holds for $\eta \leq \eta_2 \approx 0.4302$. Since $\eta_2>\eta_1$, this holds in our parameter range so $H(\Gamma)$ satisfies case (2). By the same argument as before, by induction it follows that $H^{1+2k}(\Gamma)$ contains a $v$-segment for $k \in \mathbb{N}$. Let $M = M_0 +2$, then $H^m(\gamma_u(z))$ contains a $v$-segment for all $m \geq M$. 

By an entirely analogous argument, showing that $h$-segments and their images under $H^{-1}$ must satisfy case (3)\footnote{Showing the equivalent to the $x_v(\eta)<q_1(\eta)$ bound for $H^{-1}$ requires only $\eta<\eta_2' \approx 0.4643$.}, given any $z'\in X$ we can find $N = N_0 +2$ such that $H^{-n}(\gamma_s(z'))$ contains a $h$-segment for any $n \geq N$. Since $z$ and $z'$ are arbitrary, this establishes \textbf{(MR)}.
\begin{figure}
    \centering  
    \begin{tikzpicture}[scale=2]
    \tikzmath{\e = 0.25;} 
    
    \fill[fill=gray!20] (0,{4*(1-\e)}) -- (4,0) -- ({4*(1-\e)},0) -- (0, {4*(1-\e)^2}) -- (0,{4*(1-\e)});
    
    \draw (0,0) rectangle (4-4*\e,4-4*\e);
    \draw (0,0) rectangle (4,4);
    
    \draw (0 , {4*(1-\e*\e)} ) -- ( {4*(1-\e)*\e*\e/(1-\e+\e*\e)} ,4);
    
    \draw (0,0) -- (4-4*\e, {4*(1-\e+\e*\e)}  );
    
    \draw  ( {4*(1-\e)*\e*\e/(1-\e+\e*\e)} ,0) -- (4-4*\e,4-4*\e);
    
    \node (a) at (-0.2,{0.5*4*(1-\e*\e)   +0.5*4  }) {3a'};
    \draw [very thick] (0, 4) -- (0 , {4*(1-\e*\e)} );
    
    \node (a) at (4-4*\e+0.2,{0.5*4*(1-\e+\e*\e)   +0.5*4*(1-\e)  }) {3b'};
    \draw [very thick] (4-4*\e, {4*(1-\e+\e*\e)}  )  --  (4-4*\e , {4*(1-\e)});
    
    \fill[fill=gray!40] ({4*(1-\e)^3/(1+(1-\e)^2)},{0}) -- ({4*(1-\e)^2/(1+(1-\e)^2)},{0}) -- (0,{4*(1-\e)^3/(1+2*(1-\e)^2)}) -- (0,{4*(1-\e)^4/(1+2*(1-\e)^2)}) -- ({4*(1-\e)^3/(1+(1-\e)^2)},{0});
    
    \node at ({4*(1-\e)^3/(1+(1-\e)^2)},-0.1) {$q_1$};
    \node at ({4*(1-\e)^2/(1+(1-\e)^2)},{-0.1}) {$q_2$};
    \node at (-0.1,{4*(1-\e)^3/(1+2*(1-\e)^2)}) {$q_3$};
    \node at (-0.1,{4*(1-\e)^4/(1+2*(1-\e)^2)}) {$q_4$};

    \fill[fill=gray!40] ({4*(1-\e)^3/(1+(1-\e)^2)},{4*(1-\e)^2/(1+(1-\e)^2)} ) -- ({4*(1-\e)^2/(1+(1-\e)^2)},{4*(1-\e)/(1+(1-\e)^2)}) -- ({4*(1-\e)^2/(1+2*(1-\e)^2)},{4*(1-\e +(1-\e)^3)/(1+2*(1-\e)^2)}) -- ({4*(1-\e)^3/(1+2*(1-\e)^2)},{4*((1-\e)^2 +(1-\e)^4)/(1+2*(1-\e)^2)}) -- ({4*(1-\e)^3/(1+(1-\e)^2)},{4*(1-\e)^2/(1+(1-\e)^2)} );
    
    \draw [dashed, thick] ({4*(1-\e)^2/(1+(1-\e)^2)},{0}) -- (0,{4*(1-\e)^3/(1+2*(1-\e)^2)});
    
    \draw [dashed, thick] (0,{4*(1-\e)^4/(1+2*(1-\e)^2)}) -- ({4*(1-\e)^3/(1+(1-\e)^2)},{0});
    
    \draw [dashed, thick]({4*(1-\e)^2/(1+(1-\e)^2)},{4*(1-\e)/(1+(1-\e)^2)}) -- ({4*(1-\e)^2/(1+2*(1-\e)^2)},{4*(1-\e +(1-\e)^3)/(1+2*(1-\e)^2)});
    
    \draw [dashed, thick] ({4*(1-\e)^3/(1+2*(1-\e)^2)},{4*((1-\e)^2 +(1-\e)^4)/(1+2*(1-\e)^2)}) -- ({4*(1-\e)^3/(1+(1-\e)^2)},{4*(1-\e)^2/(1+(1-\e)^2)} );
    
    \node at ({4*(1-\e)^2/(1+(1-\e)^2)},0.2) {$Q_1$};
    \node at ({4*(1-\e)^3/(1+(1-\e)^2)},{4*((1-\e)^2 +(1-\e)^4)/(1+2*(1-\e)^2) +0.1} ) {$H(Q_1)$};
    
     \node[scale=2] at (4-5*\e,2*\e) {$A_2$};
    
    \end{tikzpicture}
    
    \caption{Diagram showing that if a segment $\Gamma$ connects 3a' to 3b' through $A_3'$, then $H(\Gamma)$ must satisfy case (2).}
    \label{fig:rmip}
\end{figure}
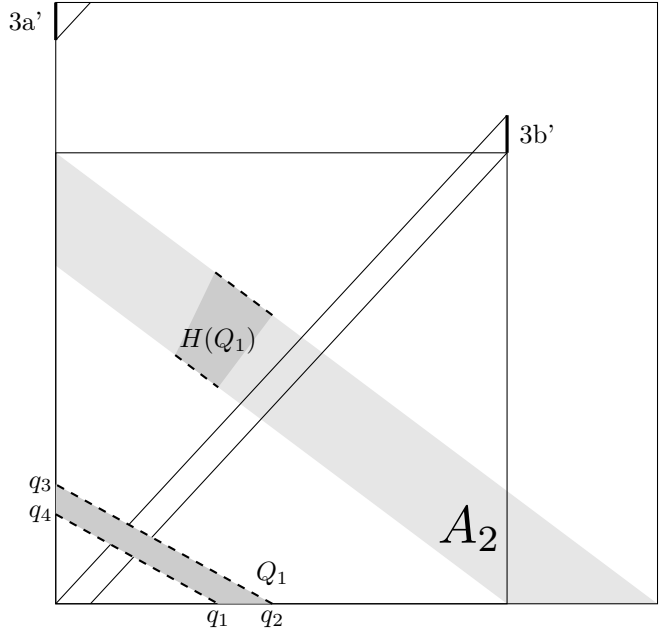
\end{proof}

We are now ready to prove the main theorem.
\begin{proof}[Proof of Theorem \ref{thm:Bernoulli}]
Noting that \textbf{(KS1)} and \textbf{(KS2)} were trivially satisfied, and the other conditions follow from Propositions \ref{prop:nonZeroChi}, \ref{prop:M}, \ref{prop:MR} over $0<\eta<\eta_1$, Theorem \ref{thm:katok-strelcyn} gives the Bernoulli property for $\eta$ over this range.
\end{proof}

%% file: sections/mixingRate.tex
The fact that we have strong expansive behaviour across the entire domain under just one iterate of $H$ allows us to deduce an exponential rate of mixing with minimal further analysis. Define the correlation function $C_n$ as in (\ref{eq:corrFunction}). We will show:

\begin{thm}
\label{thm:ourMixingRate}

Let $0<\eta<\eta_1$. There exists constants $c_1,c_2 >0$ such that $|C_n(\phi,\psi,H,\mu)|<c_1 e^{-c_2 n}$ for all H\"older continuous observables $\phi,\psi$. That is, we have exponential decay of correlations.
\end{thm}

We will follow a scheme outlined in \cite{chernov_billiards_2005}, which gives easily verifiable conditions under which a system exhibits exponentially decaying return times to a subset $\Lambda$, and subsequently exponential decay of correlations by construction of a Young Tower. We first list some basic properties for systems amenable to the scheme, paraphrased from \cite{chernov_billiards_2005}.

Let $M$ be an open domain in a 2D $C^\infty$ compact Riemannian manifold $\mathcal{M}$ with or without boundary, $f:M \rightarrow M$.
\begin{enumerate}[label={\bfseries (CZ\arabic*):}]
    \item Smoothness. The map $f$ is a $C^2$ diffeomorphism of $M\setminus \mathcal{S}$ onto $f( M \setminus \mathcal{S})$, where $\mathcal{S}$ is a closed set of zero Lebesgue measure.
    \item Hyperbolicity. At any $x\in M'\subset M$ where $Df_x$ exists, there exists two families of cones $C_x^u$ (unstable) and $C_x^s$ (stable) such that $Df_x (C_x^u) \subset C_{f(x)}^u$ and $Df_x (C_x^s) \supset C_{f(x)}^s$. There exists a constant $\lambda>1$ such that
    \[ ||Df_x(v)|| \geq \lambda||v|| \, \forall v \in C_x^u \quad \text{and} \quad||Df_x^{-1}(v)|| \geq \lambda||v|| \, \forall v \in C_x^s. \]
    These families of cones are continuous on $M'$, and the angle between $C_x^u$ and $C_x^s$ is bounded away from zero. For any $f$-invariant measure $\mu'$, at almost every $x\in M$ we have non-zero Lyapunov exponents and can define local unstable and stable manifolds $W^u(x)$, $W^s(x)$.
    \item SRB measure. The map $f$ preserves an measure $\mu$ whose conditional distributions on unstable manifolds are absolutely continuous, and is mixing.
    \item Distortion bounds. Let $\lambda(x)$ denote the factor of expansion on the unstable manifold $W^u(x)$. If $x,y$ belong to an unstable manifold $W^u$ such that $f^n$ is defined and smooth on $W^u$, then
    \[ \log \prod_{i=0}^{n-1} \frac{\lambda(f^ix)}{\lambda(f^iy)} \leq \alpha( \mathrm{dist} \left(  f^nx,f^ny \right) ) \]
    where $\alpha(\cdot)$ is some function, independent of $W^u$, with $\alpha(s) \rightarrow 0$ as $s \rightarrow 0$.
    \item Bounded Curvature. The curvature of unstable manifolds is uniformly bounded by a constant $B\geq 0$.
    \item Absolute continuity. If $W_1,W_2$ are two small unstable manifolds close to each other, then the holonomy map $h:W_1 \rightarrow W_2$ (defined by sliding along stable manifolds) is absolutely continuous with respect to the induced Lebesgue measures $\nu_{W_1}$ and $\nu_{W_2}$, and its Jacobian is bounded:
    \[ \frac{1}{C'} \leq \frac{\nu_{W_2}(h(W_1'))}{\nu_{W_1}(W_1')} \leq C' \]
    for some $C'>0$, where $W_1'\subset W_1$ is the set of points where $h$ is defined.
    \item Structure of the singularity set. For any unstable curve $W\subset M$ (a curve whose tangent vectors lie in unstable cones) the set $W \cap \mathcal{S}$ is finite or countable and...\footnote{There is an additional requirement in the countable case, irrelevant for our particular singularity set.}
\end{enumerate}
Denote the length of an line segment $W$ by $|W|$. Denote the connected components of $W\cap (M\setminus \mathcal{S})$ by $W_i$. We are now ready to give the result from \cite{chernov_billiards_2005}, specifically their Theorem 10 with $m=1$.

\begin{thm}[\citeauthor{chernov_billiards_2005}]
\label{thm:chernovZhang}

Let $f$ be defined on a 2D manifold $\mathcal{M}$ and satisfy the requirements \textbf{(CZ1-7)}. Suppose
\begin{equation}
    \label{eq:oneStepExpansion}
     \liminf_{\delta \rightarrow 0} \sup_{W: |W|<\delta} \sum_i \lambda_i^{-1} < 1     
\end{equation}
where the supremum is taken over unstable manifolds $W$ and $\lambda_i$ denotes the minimal expansion factor on $W_i$. Then the map $f:M \rightarrow M$ enjoys exponential decay of correlations. 
\end{thm}

\subsection*{Application to our maps}
We now apply this scheme directly to establish exponential decay of correlations for $H$.

\begin{proof}[Proof of Theorem \ref{thm:ourMixingRate}]
Take $M = \mathcal{M} = \tor$ and $f=H$. Starting with \textbf{(CZ1)}, take $\mathcal{S} = \mathcal{D}$ as defined in the introduction and let $M' = \tor \setminus \mathcal{D}$. Clearly $H:M' \rightarrow H(M')$ is a $C^2$ diffeomorphism and $\mu(\mathcal{D}) =0$. Moving onto \textbf{(CZ2)}, take $C_x^u = \mathcal{C}$ and $C_x^s = \mathcal{C}'$ for all $x \in M'$. Clearly these are continuous over $M'$ with cone invariance, expansion\footnote{Expansion in the $||\cdot||_2$ norm follows from a similar argument.}, transversality shown in section \ref{sec:lyp}. \textbf{(CZ3)} follows from Theorem \ref{thm:Bernoulli}, noting that Bernoulli implies strong mixing in the ergodic hierarchy. Next \textbf{(CZ4)}, \textbf{(CZ5)} follow from piecewise linearity of $H$ and \textbf{(CZ6)} follows from \textbf{(KS1-3)}. Finally since vectors tangent to $\mathcal{D}$ lie in $\mathcal{C}'$, unstable curves $W$ (with tangent vectors in $\mathcal{C}$) meet $\mathcal{D}$ transversally. Since $\mathcal{D}$ is a finite collection of segments, $W\cap \mathcal{D}$ is finite, satisfying \textbf{(CZ7)}.

It remains to show the one step expansion condition (\ref{eq:oneStepExpansion}). Note that by inspection of the partition $A_i$ in Figure \ref{fig:Apartition}, we can pick $\delta$ sufficiently small so that any unstable manifold of length $\delta$ has at most three intersections with $\mathcal{D}$, giving four connected components $W_i = W\cap{A_i}$. Note that each expansion factor $\lambda_i$ is then bounded from below by $\mathcal{K}_i$ which are defined in the same way as the $K_i$ in section \ref{sec:man}, only using the $||\cdot||_2$ norm rather than $||\cdot||_\infty$. So (\ref{eq:oneStepExpansion}) holds provided that
$\sum_{i=1}^4 \mathcal{K}_i^{-1} < 1$, which can be rewritten as
\begin{equation}
\label{eq:KinequalityMixingRates}
     \frac{1}{\mathcal{K}_1(\eta)} < 1 - \frac{1}{\mathcal{K}_2(\eta)} - \frac{1}{\mathcal{K}_3(\eta)} - \frac{1}{\mathcal{K}_4(\eta)}.
\end{equation}
Note the similarity with inequality (\ref{eq:expansionInequality}) in section \ref{sec:man}. Considering expansion under the $||\cdot||_2$ norm results in a less stringent bound on the parameter space, so that this holds over $0<\eta<\eta_1$ as required. By Theorem \ref{thm:chernovZhang}, $H$ then exhibits exponential decay of correlations over this parameter range.
\end{proof}

%% file: sections/2Dparam.tex
In this section we generalise our results to the full $0<\xi,\eta<1$ parameter space. We begin by making a small adjustment to Lemmas \ref{lemma:unstableGrowth} and \ref{lemma:stableGrowth} in the context of the $\xi=\eta$ parameter space, which will allow us to establish a larger mixing window in section \ref{sec:2DmixingProperties}.

\subsection{Weakening the growth condition}
We begin by weakening condition (\ref{eq:expansionInequality}) in section \ref{sec:man}. Let $\Gamma_p$ be a line segment which has simple intersection with each of the $A_i$. For each $M_j$ let $K_j(\eta,v) = || M_j v ||_\infty ||v||_\infty^{-1}$ be the expansion factor for the matrix $M_j$ in the direction $v\in \mathcal{C}$. Condition (\ref{eq:expansionInequality}) can be rewritten in the form
\[ \sum_{i=1}^4 \sup_{v \in \mathcal{C}} \frac{1}{K_i(\eta,v)} <1\]
which assumes $\Gamma$ to have the least expansive gradient in each $A_i$. But since the gradient of $\Gamma$ is constant, an equally valid (and weaker) condition is given by
\begin{equation}
    \label{eq:newMixingInequality}
    \sup_{v \in \mathcal{C}} \sum_{i=1}^4 \frac{1}{K_i(\eta,v)} <1.
\end{equation}
Unit vectors in $\mathcal{C}$ are of the form $(k,1)^T$ for $k_0 \leq k \leq k_1$ with $k_0 = \eta/\left(\eta-1\right)$, $k_1=1$. For each $i$ let $M_i = \begin{pmatrix} a_i & b_i \\ c_i & d_i \end{pmatrix}$, then
\begin{equation*}
    \begin{split}
        \sum_{i=1}^4 \frac{1}{K_i(\eta,v)} & = \sum_{i=1}^4 \frac{1}{|c_ik + d_i|} \\
        & = \frac{1}{c_1k + d_1} +\frac{1}{-c_2k - d_2} +\frac{1}{-c_3k - d_3} +\frac{1}{c_4k + d_4}\\
        & =: \Phi(\eta,k)
    \end{split}
\end{equation*}
where we have used the fact that $M_2$ and $M_3$ are orientation reversing. Now
\[ \frac{\partial^2 \Phi}{\partial k ^2} =  \frac{2 c_1^2}{(c_1k + d_1)^3} +\frac{2c_2^2}{(-c_2k - d_2)^3} +\frac{2c_3^2}{(-c_3k - d_3)^3} +\frac{2c_4^2}{(c_4k + d_4)^3} \]
which, by comparing with the terms of $\Phi(\eta,k)$, is clearly positive. Hence for each $\eta$, $\Phi$ as a function in $k$ is convex, giving
\[\sup_{v \in \mathcal{C}} \sum_{i=1}^4 \frac{1}{K_i(\eta,v)} = \sup_{k_0 \leq k \leq k_1} \Phi(\eta,k) = \max \{ \Phi(\eta,k_0), \Phi(\eta, k_1)\}. \]
Over $0\leq \eta <\frac{1}{2}$ we have that $\Phi(\eta,k_0) > \Phi(\eta, k_1)$ so that (\ref{eq:newMixingInequality}) (and by extension the growth lemma in forwards time) holds over $0<\eta<\eta_3$ where $\eta_3 \approx 0.2389$ solves the equation $\Phi(\eta,k_0) = 1$. The growth lemma in backwards time requires the same inequality.

\subsection{Parameter space symmetries}
Note that the system of maps $H_{(\xi,\eta)} = G_\xi \circ F_\eta$ given in the introduction is well defined and incorporates two non-monotonic shears for all $0< \xi, \eta <1$. Two symmetries exist which allow us to reduce this parameter space by a factor of four. Firstly consider $\sigma_1(\xi,\eta) = (\eta,\xi)$, reflection in the line $\eta=\xi$. We claim that
\[ S_1 \circ G_\xi \circ F_\eta = F_\xi \circ G_\eta \circ S_1 \]
where $S_1 :\tor \rightarrow \tor$ maps $(x,y) \mapsto (y,x)$. This follows from the fact that $S_1(A_j)$ = $G_\eta^{-1}(S_1(R_j))$ for $j=1,\dots,4$ and the definitions of $F$, $G$ given in the introduction. Let $\mathcal{H} = F \circ G$ (shearing vertically first instead of horizontally) then it follows that we have a semi-conjugacy between $H_{(\xi,\eta)}$ and $\mathcal{H}_{(\eta,\xi)} = \mathcal{H}_{\sigma_1(\xi,\eta)}$. Clearly $H$ and $\mathcal{H}$ share the same mixing properties, so mixing properties of $H_{\sigma_1(\xi,\eta)}$ follow from those of $H_{(\xi,\eta)}$.

Similarly take $\sigma_2(\xi,\eta) = (1-\eta,1-\xi)$, reflection in the line $\eta=1-\xi$. One can verify that
\[ S_2 \circ G_\xi \circ F_\eta = F_{1-\xi}^{-1} \circ G_{1-\eta}^{-1} \circ S_2 \]
where $S_2 :\tor \rightarrow \tor$ maps $(x,y) \mapsto (1-y,1-x)$, noting that $S_2(A_j) = G_{1-\eta}(S_2(R_j))$. This gives $H_{(\xi,\eta)}$ conjugate to $H_{\sigma_2(\xi,\eta)}^{-1}$, which has the same mixing properties as $H_{\sigma_2(\xi,\eta)}$.

Taking both of these symmetries into account, we need only study the reduced parameter space $\mathcal{P}$ defined by $\xi \leq \eta \leq 1-\xi$ with $0 <\xi \leq \frac{1}{2}$.

\subsection{Elliptic islands}
\label{sec:islands}
We state without proof a generic result on elliptic islands for piecewise linear toral automorphisms.
\begin{prop}
\label{prop:islandExists}
Let $H$ be a piecewise linear, continuous, area-preserving toral map with singularity set $\mathcal{D}$. Suppose $H$ admits an order $n$ periodic orbit $\{ z_1, z_2, \dots, z_n \}$ such that the associated cocycle $M = DH_{z_1}^n$ satisfies $|\mathrm{tr}(M)|<2$ and
$\mathrm{dist}(z_k,\mathcal{D})>0$ for $k=1,\dots,n$. Then there exists an ellipse $E$ centred at $z_1$ such that $H^n(E) = E$.
\end{prop}
We now apply the result to three periodic orbits of $H_{(\xi,\eta)}$.
\begin{cor}
$H$ exhibits elliptic islands of positive measure over the following parameter spaces:
\begin{enumerate}[label={ ($\mathcal{I}_\arabic*$):}]
    \item $ \frac{1}{2}<\eta\leq 1-\xi$ for $0 \leq \xi <\frac{1}{2}$,
    \item $0 < \xi < \min \left\{  1- \frac{1}{3\eta} , \frac{8\eta^{3}-22\eta^{2}+18\eta+\sqrt{4\eta^{3}-4\eta^{2}+1}-5}{2\left(4\eta^{3}-9\eta^{2}+7\eta-2\right)} \right\}$,
    \item $\max \left\{ \frac{1}{3-3\xi}, \frac{2\xi^{2}-4\xi+1}{2\xi^{2}-3\xi+1} \right\} < \eta < \frac{1}{2}$.
\end{enumerate}
\end{cor}
\begin{proof}
Starting with $\mathcal{I}_1$, consider the periodic orbit $\{ z_1, z_2 \}$ where
\[ z_1 = \left( \frac{- 2 \xi^{2} \eta + 5 \xi \eta - \xi - 3 \eta + 1}{4 \xi \eta - 4 \eta + 1}
, \frac{- 2  \xi \eta^{2} + 3  \xi \eta + 2  \eta^{2} - 4  \eta + 1 }{4  \xi \eta - 4  \eta + 1 }
  \right) \]
and
\[ z_2 = \left(- \frac{4 \eta \left(2 \xi^{2} - 3 \xi + 1\right)}{4 \xi \eta - 4 \eta + 1},\frac{- 2 \xi \eta^{2} + 5 \xi \eta + 2 \eta^{2} - 5 \eta + 1}{4 \xi \eta - 4 \eta + 1}\right).\]
We claim that for $(\xi,\eta) \in \mathcal{I}_1$ both $z_1 = (x_1,y_1)$ and $z_2 = (x_2,y_2)$ are contained in the interior of $A_3$, i.e. both $F(z_k)$ are in $R_3$. Now $F(x_1,y_1) = (x_2,y_1)$ and $F(x_2,y_2) = (x_1,y_2)$ so we require $0<x_k<1-\xi$ and $1-\eta<y_k<1$, which is easily verified for $(\xi,\eta) \in \mathcal{I}_1$. It follows that $\mathrm{dist}(z_k,\mathcal{D})>0$ and the associated cocycle is $M_3 M_3$.
We remark that $\mathrm{tr}(M^2) = (\mathrm{tr}M)^2 -2 \,\mathrm{det}M$ so that for area preserving matrices $M$, we have $|\mathrm{tr}(M^2)|<2 \iff |\mathrm{tr}M|<2$. Hence the conditions listed in Proposition \ref{prop:islandExists} are verified provided that $\left| 2-1/(\eta-\eta\xi) \right| < 2$, i.e. $4\eta(1-\xi)>1$, which clearly holds over $\mathcal{I}_1$.

The analysis for $\mathcal{I}_2$ and $\mathcal{I}_3$ is analogous. They correspond to islands around period 6 orbits with itinerary $A_3, A_3, A_1, A_3, A_3, A_1$. The condition on the trace of the associated cocycle gives $\xi<1-1/(3\eta)$, equivalently $\eta>1/(3-3\xi)$. The other bounds on $\mathcal{I}_2$, $\mathcal{I}_3$ come from requiring $\mathrm{dist}(z_k,\mathcal{D})>0$.
\end{proof}
The parameter regions $\mathcal{I}_n$ and their symmetries under $\sigma_1$, $\sigma_2$ are shown in Figure \ref{fig:fullParameterSpace}. These are the three largest (in terms of proportion of the parameter space) elliptic island families over $\mathcal{P}$ but do not constitute an exhaustive list. Numerical evidence suggests that the parameter space close to $\mathcal{I}_2$ and $\mathcal{I}_3$ contains parameters where $H_{(\xi,\eta)}$ is globally hyperbolic, and others where it admits other families of elliptic islands.

\subsection{Mixing properties}
\label{sec:2DmixingProperties}
In this section we generalise our approach for proving mixing properties over the line $\eta=\xi$ to subsets of $\mathcal{P}$. Inequalities on generalised expansion factors dictate where in $\mathcal{P}$ we can establish hyperbolicity, \textbf{(MR)}, and exponential decay of correlations. Starting with hyperbolicity, across $\mathcal{P}$ the traces of the $M_j(\xi,\eta)$ satisfy $|\mathrm{tr}M_j|>2$ for $j=1,2,4$. For $M_3$ we have $\mathrm{tr}M_2(\xi,\eta) = 2-1/\left(\eta-\eta\xi \right)$ which has absolute value greater than 2 provided that $1/\left(\eta-\eta\xi \right)>4$, i.e. for $\eta< 1/\left(4-4\xi\right)$. Let $\mathcal{P}'$ denote the points in $\mathcal{P}$ for which this inequality is satisfied. We remark that the cone $\mathcal{C}$ bounded by the unstable eigenvectors of $M_2$ and $M_3$, containing those of $M_1$ and $M_4$, is invariant and expanding for parameter values in $\mathcal{P}'$. The cone $\mathcal{C}'$ for $H^{-1}$ is similar, bounded by the stable eigenvectors of $M_2$ and $M_3$. Under the $||\cdot||_\infty$ norm, the cone boundaries of $\mathcal{C}$ are given by the unit vectors $(k_0,1)^T$ and $(k_1,1)^T$, where
\[ k_0(\xi,\eta) = \frac{-2\xi}{1+\sqrt{1-4\xi+4\xi\eta}}<0\quad \text{ and }  \quad k_1(\xi,\eta) = \frac{2-2\xi}{1+\sqrt{1-4\eta+4\xi\eta}}>0. \]
The cone boundaries of $\mathcal{C}'$ are given by the unit vectors $(1,m_0)^T$ and $(1,m_1)^T$, where
\[ m_0(\xi,\eta) = \frac{\sqrt{4 \xi \eta-4\xi+1}-1}{2\xi} \quad \text{ and } \quad m_1(\xi,\eta) = \frac{\sqrt{4 \xi \eta-4\eta+1}-1}{2 \xi-2}. \]

As before, write the components of $M_j$ as $a_j,\dots,d_j$ then the expansion factor $K_j(\xi,\eta,k)$ of the matrix $M_j$ in the direction $(k,1)^T \in \mathcal{C}$ is given by $|c_j k+d_j|$. Noting that each matrix has determinant 1, the expansion factor $\mathscr{K}_j(\xi,\eta,m)$ of the matrix $M_j^{-1}$ in the direction $(1,m)^T\in \mathcal{C}'$ is given by $|d_j-b_jm|$. Let \[\Phi(\xi,\eta,k) = \sum_{j=1}^4 \frac{1}{K_j(\xi,\eta,k)} \quad \text{ and } \quad \Psi(\xi,\eta,m) = \sum_{j=1}^4 \frac{1}{\mathscr{K}_j(\xi,\eta,m)},\]
then by the same reasoning as before, the growth lemma for $H$ requires $\max \{ \Phi(\xi,\eta,k_0), \Phi(\xi,\eta, k_1) \}<1$ and the growth lemma for $H^{-1}$ requires $\max \{ \Psi(\xi,\eta,m_0), \Psi(\xi,\eta, m_1) \}<1$. 

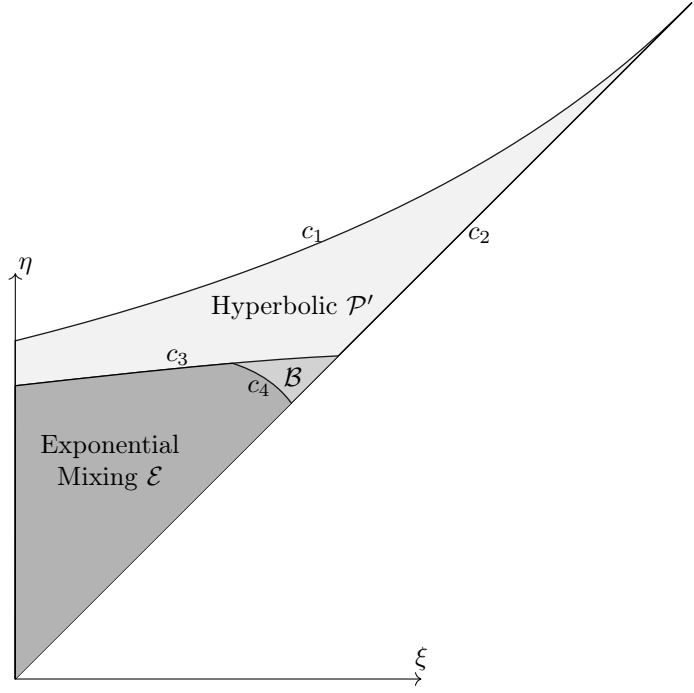
\begin{figure}
    \centering
    \begin{tikzpicture}[scale=1.8]
    \node at (2.435,2.497) {
    
    \definecolor{white}{RGB}{255,255,255}

\def \globalscale {0.231*18/25}
\begin{tikzpicture}[y=0.80pt, x=0.80pt, yscale=-\globalscale, xscale=\globalscale, inner sep=0pt, outer sep=0pt]
\begin{scope}[scale=2.000]
    \path[fill=white,rounded corners=0.0000cm] (0.0000,0.0000) rectangle (1000.0000,1000.0000);

        \path[fill=gray!10, draw=black,draw opacity=0.900,line cap=round,line join=round,line width=0.500pt,miter limit=10.00] (18.7552,499.4939) -- (18.7552,499.4939) -- (69.7492,486.3978) -- (118.8147,473.1053) -- (166.1715,459.5791) -- (211.8196,445.8400) -- (256.0030,431.8331) -- (298.7218,417.5742) -- (339.9759,403.0824) -- (380.0095,388.2874) -- (418.5784,373.2981) -- (456.1710,357.9393) -- (492.5429,342.3208) -- (527.9385,326.3500) -- (562.3576,310.0339) -- (595.8003,293.3813) -- (628.2666,276.4038) -- (660.0005,258.9780) -- (690.7580,241.2448) -- (720.7832,223.0721) -- (750.3202,204.3012) -- (779.1249,185.0771) -- (807.1972,165.4026) -- (834.7813,145.0983) -- (861.6331,124.3370) -- (887.5085,103.3301) -- (912.4074,82.1189) -- (936.3299,60.7499) -- (959.2760,39.2754) -- (980.2448,18.7552) -- (18.7308,980.2692);

        \path[draw=black,draw opacity=0.900,line cap=round,line join=round,line width=0.500pt,miter limit=10.00] (18.7308,980.2692) -- (18.7308,980.2692) -- (980.2692,18.7308);

        \path[draw=black,draw opacity=0.900,line cap=round,line join=round,line width=0.500pt,miter limit=10.00] (18.7308,980.2692) -- (18.7308,980.2692) -- (18.7308,499.5000);

        \path[fill=gray!35, draw=black,draw opacity=0.900,line cap=round,line join=round,line width=0.500pt,miter limit=10.00] (477.5273,520.9160) -- (477.5273,520.9160) -- (437.9766,522.9953) -- (394.0312,525.9383) -- (341.2969,530.1305) -- (270.9844,536.5282) -- (192.8204,544.3829) -- (92.2250,555.2125) -- (18.9408,563.4061) -- (18.7308,980.2692);

        \path[fill=gray!60,draw=black,draw opacity=0.900,line cap=round,line join=round,line width=0.500pt,miter limit=10.00] (411.0892,587.8282) -- (411.0892,587.8282) -- (405.2894,580.5270) -- (399.0232,573.6096) -- (394.0313,568.6966) -- (385.6407,561.4220) -- (372.8654,552.2344) -- (358.8750,544.1701) -- (344.9283,537.7089) -- (333.3317,533.2829) -- (327.5931,531.3604) -- (270.9844,536.5282) -- (192.8204,544.3829) -- (92.2250,555.2125) -- (18.9408,563.4061) -- (18.7308,980.2692);

\end{scope}

\end{tikzpicture}
   };

\draw[->] (0,0) -- (0,3);
\draw[->] (0,0) -- (3,0);

\node at (3.43,3.28) {$c_2$};
\node at (2.05,2.75) {Hyperbolic $\mathcal{P}'$};

\node at (2.2,3.29) {$c_1$};

\node at (1.2,2.38) {$c_3$};
\node at (1.8,2.15) {$c_4$};

\node at (2.05,2.23) {$\mathcal{B}$};

\node at (0.7,1.6) {\begin{tabular}{c}
     Exponential \\
     Mixing $\mathcal{E}$
\end{tabular}};

\node at (3,0.15) {$\xi$};
\node at (0.075,3.05) {$\eta$};

    \end{tikzpicture}
    \caption{Plot of analytical results. The curves $c_1$ and $c_2$ define $\mathcal{P}'$, $c_3$ defines $\mathcal{B}\subset \mathcal{P}'$, $c_3$ and $c_4$ define $\mathcal{E}\subset\mathcal{B}$, on which $H$ is respectively hyperbolic, mixing, and exhibits exponential decay of correlations. Note that $c_3$ meets $c_2$ at the point $(\eta_3,\eta_3)$ defined before.}
    \label{fig:final_results}
\end{figure}

Finally for each $j$ define 
\[ \mathcal{K}_j(\xi,\eta,v) = \frac{|| M_j v||_2}{||v||_2}, \]
the expansion factor of $M_j$ in the direction $v \in \mathcal{C}$ using the euclidean norm. We are now ready to state the result on mixing results over $\mathcal{P}$.
\begin{thm}
\label{thm:2dMixingResults}
\begin{itemize}
    Let $H$ be defined by parameter values $(\xi,\eta) \in \mathcal{P}$.
    \item For $(\xi,\eta) \in \mathcal{P}'$, $H$ is non-uniformly hyperbolic.
    \item For $(\xi,\eta)$ satisfying $\max \{ \Phi(k_0), \Phi(k_1), \Psi(m_0), \Psi(m_1)\}<1$, $H$ is Bernoulli.
    \item For $(\xi,\eta)$ satisfying $ \sum_j \frac{1}{\inf_{v\in \mathcal{C}} \mathcal{K}_j(\xi,\eta,v)} <1 $, $H$ exhibits exponential decay of correlations.
\end{itemize}
    The results are shown graphically in Figure \ref{fig:final_results}.
\end{thm}

\begin{proof}

The argument is similar to that given in the proofs of Theorems \ref{thm:Bernoulli} and \ref{thm:ourMixingRate}, requiring only minor adjustments. One can verify that the chain of inequalities (\ref{eq:coneInequality}) holds for all $(\xi,\eta) \in \mathcal{P}'$ so that $\mathcal{C}$ is invariant. Similarly one can verify that each of the $M_j$ expands vectors parallel to the cone boundaries, so $\mathcal{C}$ is expanding. Existence of this invariant expanding cone implies non-zero Lyapunov exponents over a full measure set, so $H$ is hyperbolic for parameter values in $\mathcal{P}'$. Moving onto proving \textbf{(M)}, Lemmas \ref{lemma:alignment}, \ref{lemma:vSegs}, and \ref{lemma:hSegs} are entirely analogous. Using the weakened condition (\ref{eq:newMixingInequality}), Lemma \ref{lemma:unstableGrowth} follows from $\max \{ \Phi(\xi,\eta,k_0), \Phi(\xi,\eta, k_1)\}<1$ and Lemma \ref{lemma:stableGrowth} follows from $\max \{ \Psi(\xi,\eta,m_0), \Psi(\xi,\eta, m_1) \}<1$. One can verify that this reduces to $\Psi(\xi,\eta, m_1)<1$, shown as the region $\mathcal{B} \subset \mathcal{P}'$ bounded by $\xi=0$, $c_2$, and the curve $c_3$ given by $\Psi(\xi,\eta, m_1)=1$ (see Figure \ref{fig:final_results}). Condition \textbf{(MR)} follows from adapting the $x_v(\eta)< q_1(\eta)$ inequality. Solving line intersection equations gives
\[ x_v(\xi,\eta) = \frac{\eta\xi(1-\xi)}{1-\eta(1-\xi)} \quad \text{and} \quad q_1(\xi,\eta) = \frac{(1-\eta)(1-\xi)^2}{1+(1-\eta)(1-\xi)} \]
so that $x_v(\xi,\eta)< q_1(\xi,\eta)$ reduces to
\[ \xi < \frac{(1-\eta)^2}{1-\eta+\eta^2} \]
which holds over $\mathcal{B}$. Again, the equivalent inequality to $x_v(\xi,\eta)< q_1(\xi,\eta)$ for $H^{-1}$ results in a less stringent condition on the parameter space, hence also holds over $\mathcal{B}$. It follows, then, that $H$ is Bernoulli over parameter values $(\xi,\eta) \in \mathcal{B}$.

Moving onto the mixing rate, \textbf{(CZ1-7)} hold by the same argument as before, noting that vectors tangent to the singularity set $\mathcal{D}$ for $H_{(\xi,\eta)}$ still lie in $\mathcal{C}'$. Similarly we can choose $\delta>0$ such that unstable manifolds $W$ of length $|W|<\delta$ have at most 3 intersections with $\mathcal{D}$, splitting $W$ into four components $W_j = W\cap A_j$. The one step expansion condition (\ref{eq:oneStepExpansion}) then follows from
\[ \sum_{j=1}^4 \sup_{v\in \mathcal{C}} \frac{1}{\mathcal{K}_j(\xi,\eta,v)} < 1,\]
i.e.
\begin{equation}
    \label{eq:curlyKinequality}
    \sum_{j=1}^4  \frac{1}{\inf_{v\in \mathcal{C}} \mathcal{K}_j(\xi,\eta,v)} < 1
\end{equation}
as required. Across $(\xi,\eta) \in \mathcal{B}$ we have that $\mathcal{K}_1(\xi,\eta,v)$ and $\mathcal{K}_2(\xi,\eta,v)$ always attain their infimum over the unstable eigenvector $v_2$ of $M_2$, $\mathcal{K}_3(\xi,\eta,v)$ and $\mathcal{K}_4(\xi,\eta,v)$ always attain their infimum over the unstable eigenvector $v_3$ of $M_3$. Hence (\ref{eq:curlyKinequality}) holds provided that $\Omega(\xi,\eta) <1$, where
\[ \Omega(\xi,\eta) = \frac{1}{\mathcal{K}_1(\xi,\eta,v_2)} + \frac{1}{\mathcal{K}_2(\xi,\eta,v_2)} + \frac{1}{\mathcal{K}_1(\xi,\eta,v_3)} + \frac{1}{\mathcal{K}_1(\xi,\eta,v_3)}.\]
Figure \ref{fig:final_results} shows the curve $c_4$ given by $\Omega(\xi,\eta) =1$ in $\mathcal{B}$, which together with $c_3$, $c_2$, $\xi=0$ give the exponential mixing window $\mathcal{E} \subset \mathcal{B}$.

\end{proof}

%% file: sections/OTM.tex
Let $H$ denote $H_{(\frac{1}{2},\frac{1}{2})}$, the map on the cusp of the hyperbolic parameter space $\mathcal{P}'$. As the composition of two orthogonal `tent' shaped shears, we will colloquially refer to this as the \emph{Orthogonal Tents Map} (OTM). It is the unique map in the full $0<\xi,\eta <1$ parameter space which is not conjugate to another $H_{(\xi,\eta)}$ and has all integer valued derivative matrices. It is also the natural extension of \citeauthor{cerbelli_continuous_2005}'s Map with two non-monotonic shears, so proving its observed hyperbolic and mixing properties is desirable, in line with other generalisations (\citealp{demers_family_2009}). We will prove the first of these, then comment on the challenges of proving the second in section \ref{sec:OTMmixing}.
\subsection{Hyperbolicity}
\label{sec:OTMhyp}
\begin{prop}
\label{prop:tentMap}
$H$ is non-uniformly hyperbolic.
\end{prop}

Let $M_j$ denote the derivative matrix $DH$ on $A_j$. These are given by
\[ M_1 = \begin{pmatrix} 1 & 2 \\ 2 & 5 \end{pmatrix}, \, M_2 = \begin{pmatrix} 1 & 2 \\ -2 & -3 \end{pmatrix}, \, M_3 = \begin{pmatrix} 1 & -2 \\ 2 & -3 \end{pmatrix}, \text{ and } M_4 = \begin{pmatrix} 1 & -2 \\ -2 & 5 \end{pmatrix}.\]
For any $z\in X'$ with $n$-step itinerary
\[  A_{j_1}, A_{j_2}, A_{j_3}, \dots, A_{j_n},\] the cocycle $DH_z^n$ is given by
\[ DH_z^n = M_{j_n} \dots M_{j_3} M_{j_2}  M_{j_1}\]
with each $j_k \in \{ 1,2,3,4 \}$. Our aim is to decompose any cocycle into hyperbolic matrices which share an invariant expanding cone.
Note that while $M_1$ and $M_4$ are hyperbolic, $M_2$ and $M_3$ are not. Hence when $M_2$ or $M_3$ appear in a cocycle at $M_{j_k}$, we must combine them with its neighbouring matrices $M_{j_{k+l}},\dots,M_{j_{k+2}}, M_{j_{k+1}}$ for some $l \in \mathbb{N}$.

Let $\mathcal{M}$ denote the countable family of matrices $\{ M_1, M_4,  M_1M_2^n, M_3M_2^n, M_4M_2^n, M_1M_3^n, M_2M_3^n, M_4M_3^n\}$ with $n\in\mathbb{N}$. We claim the following:
\begin{lemma}
\label{lemma:itineraries}
At almost every $z$, the cocycle $DH_z^n$ can be decomposed into blocks from $\mathcal{M}$.
\end{lemma}

\begin{lemma}
\label{lemma:cone}
The matrices in $\mathcal{M}$ admit an invariant expanding cone $\mathcal{C}$.
\end{lemma}

Proposition \ref{prop:tentMap} follows from the two lemmas. At any $z$ satisfying Lemma \ref{lemma:itineraries}, by Lemma \ref{lemma:cone} we can take any $v_0 \in \mathcal{C}$ to achieve exponential growth of $|| DH_z^n v_0 ||$ with $n$. We will prove Lemma \ref{lemma:itineraries} here, the proof of Lemma \ref{lemma:cone} can be found in the appendix.

\begin{proof}[Proof of Lemma \ref{lemma:itineraries}]

It is sufficient to show that itineraries cannot get \emph{trapped} in $A_2$ or $A_3$, barring some set of zero measure. We will consider the set $A_3$, with the argument for $A_2$ being entirely analogous. In particular we will show that $\mu(B_n) \rightarrow 0$ as $n\rightarrow \infty$ where $B_n = \{z' \in A_3 \,| \, H^k(z') \in A_3 \text{ for all } 1 \leq k \leq n \}$.

Let $\mathcal{H} = F \circ G$. For any $z'\in A_3$,
\begin{equation*}
\begin{split}
    H^k(z') \in A_3 \text{ for all } 1 \leq k \leq n & \iff (G \circ F)^k(z') \in A_3 \text{ for all } 1 \leq k \leq n \\
    & \iff [F \circ (G \circ F)^k ] (z') \in R_3 \text{ for all } 1 \leq k \leq n \\
    & \iff [(F \circ G)^k \circ F](z') \in R_3 \text{ for all } 1 \leq k \leq n \\
    & \iff \mathcal{H}^k(z) \in R_3 \text{ for all } 1 \leq k \leq n
\end{split}
\end{equation*}
where $z= F(z') \in R_3$. Hence recurrence in $A_3$ under $H$ can be understood by instead studying recurrence in $R_3$ under $\mathcal{H}$. Letting $\mathcal{B}_n =\{z \in R_3 \,| \, \mathcal{H}^k(z) \in R_3 \text{ for all } 1 \leq k \leq n \}$, by the above we have $\mathcal{B}_n = F(B_n)$ and $\mu(B_n) = \mu(\mathcal{B}_n)$ since $F$ preserves $\mu$. The simpler geometry of $R_3$ makes this a convenient choice. Iteratively define $U_1 = \mathcal{H}(R_3) \cap R_3$, $U_n = \mathcal{H}(U_{n-1}) \cap R_3 $ so that $\mathcal{B}_n = \mathcal{H}^{-n}(U_n)$. Since $\mathcal{H}$ preserves $\mu$, we have $\mu(\mathcal{B}_n) = \mu(U_n)$. Let $V = \mathcal{H}^{-1}(R_3) \cap R_3$ be the set of points in $R_3$ which stay in $R_3$. An equivalent definition for the $U_n$ is $U_1 = \mathcal{H}(V)$, $U_n = \mathcal{H}(U_{n-1} \cap V)$. Restricting to $V$ in this way is beneficial as $\mathcal{H}|_{V} : V \rightarrow R_3$ is an affine transformation, mapping quadrilaterals to quadrilaterals. The sets $V = V_1 \cup V_2$ and $U_1 = P_1 \cup Q_1$ are shown in Figure \ref{fig:UV}, both composed of two quadrilaterals with corners on $\partial R_3$. Note that $V_1$, $P_1$ share the corners $p_1^1 = \left(1/4,1/2 \right)$, $p_1^3 = \left(0,3/4 \right)$ and $V_2$, $Q_1$ share the corners $q_1^1 =  \left(1/4, 1\right)$, $q_1^3 = \left(1/2,3/4 \right)$, all of which are periodic with period 2.

\begin{figure}
    \centering
    \begin{tikzpicture}
    
    \node[scale=1.2] at (-4,0) {
    \begin{tikzpicture}
    
    \draw (0,5) rectangle (5,10);
    \filldraw[pattern = north west lines] (0,10) -- (2.5,5) -- (10/6,5) -- (0,7.5) -- (0,10);
    \filldraw[pattern = north west lines] (5,5) -- (2.5,10) -- (5-10/6,10) -- (5,7.5) -- (5,5);
    
    \filldraw[fill = gray!80, opacity=0.5] (0,10) -- (2.5,10) -- (5,50/6) -- (5,7.5) -- (0,10);
    
    \filldraw[fill = gray!80, opacity=0.5] (5,5) -- (2.5,5) -- (0,20/3) -- (0,7.5) -- (5,5);
        
        
    \node[scale = 1.2] at (4.7,9.7) {$R_3$};    
    
    \node at (1,6) {$\bullet$};
    \node at (0.8,5.8) {$r_1$};
    
    \node  at (10/6,20/3) {$\bullet$};
    \node  at (10/6+0.25,20/3+0.25) {$r_1'$};
    
     \node at (10/3,50/6) {$\bullet$};
     \node at (10/3-0.2,50/6-0.2) {$s_1$};
    
     \node  at (20/5,9) {$\bullet$};
   \node  at (20/5+0.2,9+0.2) {$s_1'$};
    
    \node at (3,5.5) {$P_1$};
    \node at (2,9.5) {$Q_1$};
    \fill[fill=white] (0.5,8) circle (6.5pt);
    \node at (0.5,8) {$V_1$};
    
    \fill[fill=white] (4.5,7) circle (6.5pt);
    \node at (4.5,7) {$V_2$};
    
    \node[scale=0.8] at (10/6,4.75) {$\frac{1}{6}$};
    \node[scale=0.8] at (2.5,4.75) {$\frac{1}{4}$};
    \node[scale=0.8] at (10/3,10.25) {$\frac{1}{3}$};
    \node[scale=0.8] at (2.5,10.25) {$\frac{1}{4}$};
    \node[scale=0.8] at (-0.15,7.5) {$\frac{3}{4}$};
    \node[scale=0.8] at (-0.15,20/3) {$\frac{2}{3}$};
    \node[scale=0.8] at (5.15,7.5) {$\frac{3}{4}$};
    \node[scale=0.8] at (5.15,50/6) {$\frac{5}{6}$};

    \end{tikzpicture}};
    
    \draw [->, thick] (-0.5,0) -- (0.5,0);
    \node[scale=1.3] at (0,0.3) {$\mathcal{H}$};
    
    \node[scale=1.2] at (4.3,0) {
    \begin{tikzpicture}
    
    \draw (0,5) rectangle (5,10);
        
    \filldraw[pattern = north west lines] (5.0, 8.0) -- (5.0, 7.5) -- (5/3, 10) -- (2.5, 10) -- (5.0, 8.0);
    \filldraw[fill = gray!80, opacity=0.5] (5.0, 8.0) -- (5.0, 7.5) -- (5/3, 10) -- (2.5, 10) -- (5.0, 8.0);
    
    \filldraw[pattern = north west lines] (0,7.5) -- (10/3,5) -- (2.5,5) -- (0,7) -- (0,7.5);
    \filldraw[fill = gray!80, opacity=0.5] (0,7.5) -- (10/3,5) -- (2.5,5) -- (0,7) -- (0,7.5);

    \draw[dashed] (0,10) -- (2.5,5) -- (10/6,5) -- (0,7.5) -- (0,10);
    \draw[dashed] (5,5) -- (2.5,10) -- (5-10/6,10) -- (5,7.5) -- (5,5);    
        
    \node[scale = 1.2] at (4.7,9.7) {$R_3$};
    
    \node[scale=0.8] at (2.5,4.8) {$p_2^1$};
    \node[scale=0.8] at (-0.2,7) {$p_2^2$};
    \node[scale=0.8] at (-0.2,7.5) {$p_2^3$};
    \node[scale=0.8] at (10/3,4.8) {$p_2^4$};
    
    \node[scale=0.8] at (2.5,10.2) {$q_2^1$};
    \node[scale=0.8] at (5.2,8) {$q_2^2$};
    \node[scale=0.8] at (5.2,7.5) {$q_2^3$};
    \node[scale=0.8] at (5-10/3,10.2) {$q_2^4$};
    
    \node at (10/14,90/14) {$\bullet$};
    \node at (10/14-0.2,90/14-0.2) {$r_2$};
    
    \node at (2,6) {$\bullet$};
    \node at (2+0.2,6+0.2) {$r_2'$};
    
    \node at (60/14,120/14) {$\bullet$};
    \node at (60/14+0.2,120/14+0.2) {$s_n$};
    
    \node at (3,9) {$\bullet$};
    \node at (3-0.2,9-0.2) {$s_n'$};

    \end{tikzpicture}};
    \end{tikzpicture}
    \caption{Left: Two subsets $V$ (patterned) and $U_1$ (grey) of $R_3$, each composed of two quadrilaterals. Right: The image $U_2 = \mathcal{H}(U_1 \cap V)$ in $R_3$, the dashed lines show the boundary of $V$.}
    \label{fig:UV}
\end{figure}
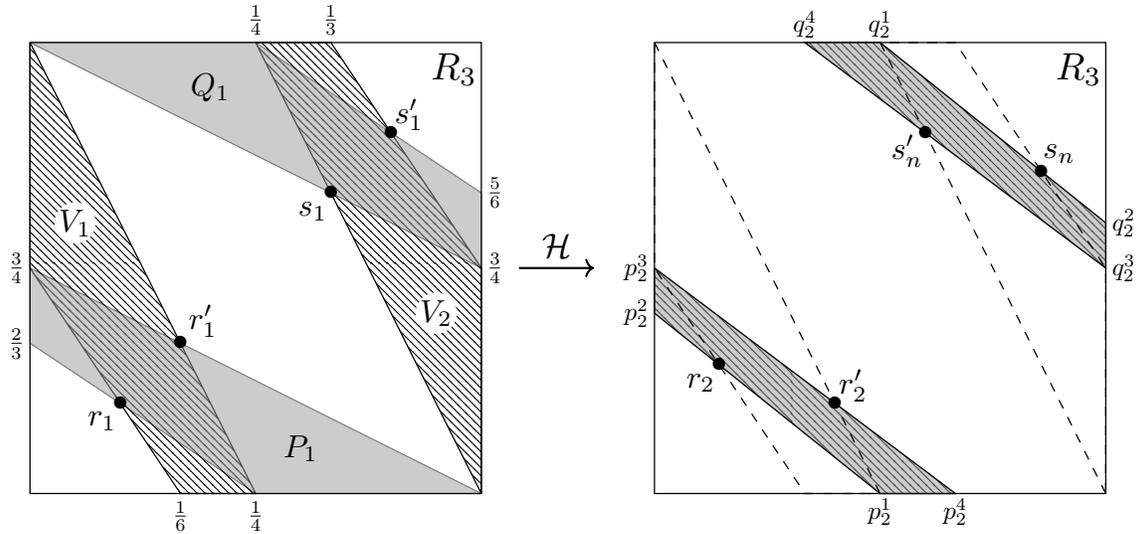

The intersection $U_1 \cap V$ is made up of two quadrilaterals $P_1 \cap V_1$ and $Q_1 \cap V_2$ with corners on the period 2 points and the points $r_1 = \left( 1/10 ,3/5 \right)$, $r_1' = \left( 1/6,2/3\right)$, $s_1 = \left( 1/3 ,5/6 \right)$, and $s_1' = \left( 2/5 ,9/10 \right)$. Mapping these quadrilaterals forward under $\mathcal{H}$ gives $U_2 = P_2 \cup Q_2$ where $P_2 = \mathcal{H}(Q_1 \cap V_2)$ and $Q_2 = \mathcal{H}( P_1 \cap V_1)$. Label the corners of these quadrilaterals by $p_2^i$ and $q_2^i$, $i=1,2,3,4$, as shown in Figure \ref{fig:UV}.

We claim that for general $n\in \mathbb{N}$, $U_n$ is made up of two quadrilaterals $P_n$, $Q_n$ with corners
\[ p_n^1 = \left( \frac{1}{4}  , \frac{1}{2} \right), \quad p_n^2 = \left( 0  , \frac{3n+1}{4n+2} \right), \quad p_n^3 = \left( 0 , \frac{3}{4} \right), \quad p_n^4 = \left( \frac{n}{4n-2}  , \frac{1}{2} \right), \]
\[ q_n^1 = \left( \frac{1}{4}  , 1 \right), \quad q_n^2 = \left( \frac{1}{2}  , \frac{3n+2}{4n+2} \right), \quad q_n^3 = \left( \frac{1}{2}  , \frac{3}{4} \right), \quad q_n^4 = \left( \frac{n-1}{4n-2}  , 1 \right), \]
labelled in the same way as the case $n=2$. $P_n \cap V_1$ will be a quadrilateral with corners $p_n^1,r_n,p_n^3,r_n'$, and $Q_n \cap V_2$ will be a quadrilateral with corners $q_n^1, s_n, q_n^3,s_n'$, where
\[ r_n = \left( \frac{1}{4n+6} , \frac{3n+3}{4n+6} \right), \quad r_n' = \left( \frac{n}{4n+2} , \frac{2n+2}{4n+2} \right), \quad s_n = \left( \frac{2n+2}{4n+6} , \frac{3n+6}{4n+6} \right), \quad s_n' = \left( \frac{n+1}{4n+2} , \frac{4n+1}{4n+2} \right) \]
can be obtained by solving the line intersection equations. One can verify that $\mathcal{H} (p_n^1) = q_{n+1}^1 $, $\mathcal{H} (r_n) = q_{n+1}^2 $, $\mathcal{H} (p_n^3) = q_{n+1}^3 $, $\mathcal{H} (r_n') = q_{n+1}^4 $, and $\mathcal{H} (q_n^1) = p_{n+1}^1 $, $\mathcal{H} (s_n) = p_{n+1}^2 $, $\mathcal{H} (q_n^3) = p_{n+1}^3 $, $\mathcal{H} (s_n') = p_{n+1}^4 $, so that $\mathcal{H} (P_n \cap V_1) = Q_{n+1}$ and $\mathcal{H} (Q_n \cap V_2) = P_{n+1}$. Hence
\begin{equation*}
    \begin{split}
        \mathcal{H}(U_n \cap V) & = \mathcal{H}\left( (P_n \cap V_1) \cup (Q_n \cap V_n) \right) \\
        & = \mathcal{H} (P_n \cap V_1) \cup \mathcal{H} (Q_n \cap V_2) \\
        & = Q_{n+1} \cup P_{n+1} \\
        & = U_{n+1}
    \end{split}
\end{equation*}
and the claim follows by induction. Now in the limit $n \rightarrow \infty$, $P_n$ limits onto the line segment joining $\left( 0, 3/4 \right)$ to $\left(1/4,1/2\right)$ and $Q_n$ limits onto the line segment joining $\left( 1/4, 1 \right)$ to $\left( 1/2, 3/4 \right)$. This give $\mu(U_n) \rightarrow 0$ as required. The preimages of these segments under $F$ are visible as the darker regions of the Poincar\'e section given in Figure \ref{fig:otmPoin}. If an orbit (like that shown in the figure) maps near to the segments, it can take arbitrarily long to escape. This gives a non-uniform spatial density for the orbit in the \emph{finite time} picture of the dynamics that a Poincar\'e section provides.

\end{proof}

\subsection{Mixing properties}
\label{sec:OTMmixing}
This approach of identifying non-hyperbolic regions $A_2$, $A_3$ and proving that itineraries cannot get trapped there is similar to the method used to prove hyperbolicity and the mixing property in \cite{myers_hill_continuous_2021}. Unlike the maps studied there, which had finite escape times from the non-hyperbolic region, here we can find positive Lebesgue measure sets which take arbitrarily long to escape. This complicates establishing the growth lemma, requiring analysis of a countably infinite partition of returns, so that the proof of the mixing property is more involved. This is the subject of current work.

Despite this, the above analysis allows us to comment on the potential mixing rate of the OTM. By the shoelace formula we can calculate $\mu(P_n) = \mu(Q_n) = n/\left(32n^2-8\right)$ so that $\mu(B_n)$, the measure of the unmixed region in $A_3$, is given by $\mu(B_n) = \mu(U_n) =  n/\left(16n^2-4\right)$. This suggests that the mixing rate is at most \emph{polynomial}, in contrast to the exponential mixing rate seen elsewhere in our parameter space. Numerical evidence supports this, and suggests exponential correlation decay rate across $\mathcal{P}'$ and the curve $c_1$ left of the OTM.

%% file: sections/discussion.tex
\subsection{Improving the (exponential) mixing windows}
Numerical results suggest mixing behaviour across all of $\mathcal{P}'$ and some way beyond. The key issue limiting our analysis from establishing mixing results over the larger parameter space is the weak hyperbolicity of $M_3$ near $\eta=1/\left(4-4\xi\right)$ and non-hyperbolicity for $\eta>1/\left(4-4\xi\right)$. There are methods for getting around this weak expansion, considering expansion over $n$ iterates and using the precise geometry of the singularity set for $H^n$ to derive stronger bounds on the growth of local manifolds. Several factors prevent the easy application of this method. Firstly, neighbouring partition elements defined by the singularity set for $H^n$ will always have inverse orientation preserving/reversing properties. This was not the case in \cite{myers_hill_continuous_2021} and was key in establishing an analogous growth lemma for piecewise linear curves rather than line segments. Secondly, considering $H^n$ with two non-monotonic shears involves working with a very complicated singularity set with $4^n$ partition elements. This, together with a two-dimensional parameter space, makes any analysis significantly more challenging.

Recall that the one step expansion condition in Theorem \ref{thm:chernovZhang} (\citealp{chernov_billiards_2005}) was the key constraint on our exponential mixing window $\mathcal{E}$. In subsequent publications this condition has been weakened, employing image coupling methods rather than construction of a Young tower. Using similar notation to Theorem \ref{thm:chernovZhang} above, the weakened condition (from \citealp{chernov_statistical_2009}) is given as follows for our map $H$. Let $W$ be an unstable curve, $W_i$ be the restriction to $A_i$, and $V_i = H(W_i)$.
The one-step expansion condition is satisfied provided that there exists $q\in (0,1]$ such that
\begin{equation}
    \label{eq:newOneStepExpansion}
     \liminf_{\delta \rightarrow 0} \sup_{W: |W|<\delta} \sum_i \left( \frac{|W|}{|V_i|} \right)^q \cdot \frac{|W_i|}{|W|} < 1,
\end{equation}
where the supremum is taken over all unstable curves $W$. This is difficult to implement for our maps as finding the precise proportions of the curve(s) that attain this supremum in each of the four partition elements $A_i$ is challenging. Three proportion tuning parameters are required, which together with $q$ and the two dimensional parameter space results in a non-trivial optimisation problem.

\subsection{Comparison with linked twist maps}

Let $F$, $G$ be the non-monotonic shears given in section \ref{sec:intro}. Taking $\tilde{F} = F$, $\tilde{G} = G$ and imposing $\tilde{F}|_{\{y>1-\eta\}} = \tilde{G}|_{\{x>1-\xi\}} = \mathrm{Id}$ gives a class of linked twist maps $\tilde{H}=\tilde{G} \circ \tilde{F}$ with known mixing properties over the $(\xi,\eta)$ parameter space (\citealp{wojtkowski_linked_1980}). Indeed, the mixing rate for $\tilde{H}$ is known to be polynomial, see \cite{sturman_rate_2013} and \cite{springham_polynomial_2014}, in contrast to the exponential rate shown seen over $\mathcal{E}$ for $H$. It is clear, then, that the shears in the annuli $\{y>1-\eta\}$, $\{x>1-\xi\}$ have a significant positive impact on this aspect of the dynamics. One might ask whether including these shears improves mixing in more general linked twist maps, for example the counter-rotating LTM $\tilde{H}^-=\tilde{G}^{-1} \circ \tilde{F}$. Letting $H^-=G^{-1}\circ F$, one can show that mapping by $H^-$ rather than $\tilde{H}^-$ does not mitigate the growth of elliptic islands. For example, over $\xi=\eta < 1/2$ the maps $H^-$ and $\tilde{H}^-$ share the same pair of elliptic islands associated with the period 2 orbit
\[ \left( \frac{(1-\eta)^2}{3-2\eta} , \frac{(1-\eta)^2}{3-2\eta}  \right) \longleftrightarrow \left( \frac{(1-\eta)(2-\eta)}{3-2\eta} , \frac{(1-\eta)(2-\eta)}{3-2\eta}  \right) \]
and the non-hyperbolic matrix
\[ M_1^- = \begin{pmatrix} 1 & \frac{1}{1-\eta} \\ \frac{-1}{1-\eta} & 1+\frac{1}{(1-\eta)^2} \end{pmatrix}. \]

%% file: sections/appendix.tex
\begin{proof}[Proof of Lemma \ref{lemma:cone}]

Parameterise the tangent space by $(v_1,v_2) \in \mathbb{R}^2$. Define $\mathcal{C}$ as the cone contained within the region $|v_2| \geq |v_1|$, bounded by and including the unstable eigenvectors of $M_4M_2$ and $M_1M_3$. As unit vectors in the $|| \cdot||_\infty$ norm, these are $v^{-} = \left(-\alpha,1\right)^T$ and $ v^+ = \left(\alpha,1\right)^T$ respectively where $\alpha = \frac{1}{2}(\sqrt{5}-1)$. We will show hyperbolicity, cone invariance, and finally norm expansion of vectors in $\mathcal{C}$ under matrices from $\mathcal{M}$.

\begin{table}[h]
    \centering
    \begin{tabular}{c|c|c|c}
        $M$ &  $(\mathrm{tr}(M))^2-4$ & $g_u(M)$ & $g_s(M)$ \\ \hline
        $M_1$ & $ 32 $ & \( \displaystyle  1 + \sqrt{2} \) & \( \displaystyle 1 - \sqrt{2} \)  \\
$M_4$ & $ 32 $ & \( \displaystyle  - \sqrt{2} - 1 \) & \( \displaystyle -1 + \sqrt{2} \)  \\
$M_1M_2^n$ & $ 4 \left(4 n + 3\right)^{2} - 4 $ & \( \displaystyle  \frac{n + \sqrt{4 n^{2} + 6 n + 2} + 1}{n + 1} \) & \( \displaystyle \frac{n - \sqrt{4 n^{2} + 6 n + 2} + 1}{n + 1} \)  \\
$M_1M_3^n$ & $ 4 \left(4 n + 3\right)^{2} - 4 $ & \( \displaystyle  \frac{5 n + \sqrt{4 n^{2} + 6 n + 2} + 1}{3 n + 1} \) & \( \displaystyle \frac{5 n - \sqrt{4 n^{2} + 6 n + 2} + 1}{3 n + 1} \)  \\
$M_2M_3^n$ & $ 64 n \left(4 n + 1\right) $ & \( \displaystyle  - \frac{n + 2 \sqrt{n \left(4 n + 1\right)} + 1}{3 n + 1} \) & \( \displaystyle \frac{- n + 2 \sqrt{n \left(4 n + 1\right)} - 1}{3 n + 1} \)  \\
$M_3M_2^n$ & $ 64 n \left(4 n + 1\right) $ & \( \displaystyle  \frac{n + 2 \sqrt{n \left(4 n + 1\right)} + 1}{3 n + 1} \) & \( \displaystyle \frac{n - 2 \sqrt{n \left(4 n + 1\right)} + 1}{3 n + 1} \)  \\
$M_4M_2^n$ & $ 4 \left(4 n + 3\right)^{2} - 4 $ & \( \displaystyle  - \frac{5 n + \sqrt{4 n^{2} + 6 n + 2} + 1}{3 n + 1} \) & \( \displaystyle \frac{- 5 n + \sqrt{4 n^{2} + 6 n + 2} - 1}{3 n + 1} \)  \\
$M_4M_3^n$ & $ 4 \left(4 n + 3\right)^{2} - 4 $ & \( \displaystyle  - \frac{n + \sqrt{4 n^{2} + 6 n + 2} + 1}{n + 1} \) & \( \displaystyle \frac{- n + \sqrt{4 n^{2} + 6 n + 2} - 1}{n + 1} \)  \\

    \end{tabular}
    \caption{Information necessary for establishing hyperbolicity of each $M \in \mathcal{M}$ and for showing that they admit an invariant cone $\mathcal{C}$.}
    \label{tab:tab1}
\end{table}

Starting with hyperbolicity, a matrix $M\in \mathcal{M}$ is hyperbolic if its trace satisfies $(\mathrm{tr}(M))^2>4$. Table \ref{tab:tab1} shows $(\mathrm{tr}(M))^2-4$ for each of the matrices, one can verify that all are positive. Hence each of the matrices $M$ have distinct unstable and stable eigenvectors, write their gradients as $g_u(M)$ and $g_s(M)$ respectively. The gradients of the cone boundaries $v^{\pm}$ are $\pm 1/\alpha$, so we have cone invariance $Mv \in \mathcal{C}$ for all $v \in \mathcal{C}$ if $|g_u(M)| \geq 1/\alpha$ and $|g_s(M)| < 1/\alpha$. Again, using Table \ref{tab:tab1}, this is easily verified. 

By cone invariance, for any $M \in \mathcal{M}$, $v = (v_1,v_2)^T \in \mathcal{C}$, the vector $(v_1',v_2')^T = M (v_1,v_2)^T$ will satisfy $|| (v_1',v_2')^T  || = |v_2'|$. Write the components of $M$ as $\left(\begin{matrix} m_1 & m_2 \\ m_3 & m_4 \end{matrix}\right)$, then the expansion in norm of unit vectors $v \in \mathcal{C}$ under $M$ is given by $|m_3 v_1 + m_4|$. The expansion factors of hyperbolic matrices over vectors in an invariant cone are always minimal on one of the cone boundaries, so the minimum expansion factor for $M$ over $\mathcal{C}$ is 
\[\min \{ |\pm \alpha m_3 + m_4| \} = \begin{dcases} | -\alpha m_3 + m_4| & \text{if } \mathrm{sgn}(m_3) = \mathrm{sgn}(m_4) \\ | \alpha m_3 + m_4| & \text{if }\mathrm{sgn}(m_3) \neq \mathrm{sgn}(m_4)
\end{dcases}   \]
since $\alpha>0$. Write this minimum expansion factor as $K(M)$. Table \ref{tab:tab2} shows the minimum expansion factors $K(M)$ for each $M \in \mathcal{M}$. All are greater than 1 so that the cone $\mathcal{C}$ is expanding as required.

\begin{table}[ht]
    \centering
    \begin{tabular}{c|c|c|c}
        $M$ & Components & $K(M)$ & $\inf_n K(M)$ \\ \hline
$M_1$ & $\left(\begin{matrix}1 & 2\\2 & 5\end{matrix}\right) $ & \( \displaystyle  6 - \sqrt{5} \) & 3.763  \\
$M_4$ & $ \left(\begin{matrix}1 & -2\\-2 & 5\end{matrix}\right) $ & \( \displaystyle  6 - \sqrt{5} \) & 3.763  \\
$M_1M_2^n$ & $(-1)^n \left(\begin{matrix}2 n + 1 & 2 n + 2\\6 n + 2 & 6 n + 5\end{matrix}\right) $ & \( \displaystyle  6 n + \left(1 - \sqrt{5}\right) \left(3 n + 1\right) + 5 \) & 6.055  \\
$M_1M_3^n$ & $(-1)^n \left(\begin{matrix}1 - 6 n & 6 n + 2\\2 - 14 n & 14 n + 5\end{matrix}\right) $ & \( \displaystyle  14 n + \left(1 - \sqrt{5}\right) \left(7 n - 1\right) + 5 \) & 11.58  \\
$M_2M_3^n$ & $(-1)^n \left(\begin{matrix}1 - 6 n & 6 n + 2\\10 n - 2 & - 10 n - 3\end{matrix}\right) $ & \( \displaystyle  10 n + \left(1 - \sqrt{5}\right) \left(5 n - 1\right) + 3 \) & 8.055  \\
$M_3M_2^n$ & $(-1)^n \left(\begin{matrix}1 - 6 n & - 6 n - 2\\2 - 10 n & - 10 n - 3\end{matrix}\right) $ & \( \displaystyle  10 n + \left(1 - \sqrt{5}\right) \left(5 n - 1\right) + 3 \) & 8.055  \\ 
$M_4M_2^n$ & $(-1)^n \left(\begin{matrix}1 - 6 n & - 6 n - 2\\14 n - 2 & 14 n + 5\end{matrix}\right) $ & \( \displaystyle  14 n + \left(1 - \sqrt{5}\right) \left(7 n - 1\right) + 5 \) & 11.58  \\
$M_4M_3^n$ & $(-1)^n \left(\begin{matrix}2 n + 1 & - 2 n - 2\\- 6 n - 2 & 6 n + 5\end{matrix}\right) $ & \( \displaystyle  6 n + \left(1 - \sqrt{5}\right) \left(3 n + 1\right) + 5 \) & 6.055  \\
    \end{tabular}
    \caption{Minimum expansion factors for each $M \in \mathcal{M}$ over the cone $\mathcal{C}$.}
    \label{tab:tab2}
\end{table}

\end{proof}

%% file: main.bbl
\begin{thebibliography}{}

\bibitem[Arnold and Avez, 1968]{arnold_ergodic_1968}
Arnold, V.~I. and Avez, A. (1968).
\newblock {\em Ergodic {Problems} of {Classical} {Mechanics}}.
\newblock Addison-Wesley.

\bibitem[Beigie et~al., 1994]{beigie_invariant_1994}
Beigie, D., Leonard, A., and Wiggins, S. (1994).
\newblock Invariant manifold templates for chaotic advection.
\newblock {\em Chaos, Solitons \& Fractals}, 4(6):749--868.

\bibitem[Burton and Easton, 1980]{burton_ergodicity_1980}
Burton, R. and Easton, R.~W. (1980).
\newblock Ergodicity of linked twist maps.
\newblock In Nitecki, Z. and Robinson, C., editors, {\em Global {Theory} of
  {Dynamical} {Systems}}, Lecture {Notes} in {Mathematics}, pages 35--49,
  Berlin, Heidelberg. Springer.

\bibitem[Cerbelli and Giona, 2005]{cerbelli_continuous_2005}
Cerbelli, S. and Giona, M. (2005).
\newblock A {Continuous} {Archetype} of {Nonuniform} {Chaos} in
  {Area}-{Preserving} {Dynamical} {Systems}.
\newblock {\em Journal of Nonlinear Science}, 15(6):387--421.

\bibitem[Cerbelli and Giona, 2008]{cerbelli_characterization_2008}
Cerbelli, S. and Giona, M. (2008).
\newblock Characterization of nonuniform chaos in area-preserving nonlinear
  maps through a continuous archetype.
\newblock {\em Chaos, Solitons \& Fractals}, 35(1):13--37.

\bibitem[Chernov, 1999]{chernov_decay_1999}
Chernov, N. (1999).
\newblock Decay of {Correlations} and {Dispersing} {Billiards}.
\newblock {\em Journal of Statistical Physics}, 94(3):513--556.

\bibitem[Chernov and Young, 2000]{chernov_decay_2000}
Chernov, N. and Young, L.~S. (2000).
\newblock Decay of {Correlations} for {Lorentz} {Gases} and {Hard} {Balls}.
\newblock In {\em Hard {Ball} {Systems} and the {Lorentz} {Gas}}, Encyclopaedia
  of {Mathematical} {Sciences}, pages 89--120. Springer, Berlin, Heidelberg.

\bibitem[Chernov and Zhang, 2005]{chernov_billiards_2005}
Chernov, N. and Zhang, H.-K. (2005).
\newblock Billiards with polynomial mixing rates.
\newblock {\em Nonlinearity}, 18:1527.

\bibitem[Chernov and Zhang, 2009]{chernov_statistical_2009}
Chernov, N. and Zhang, H.-K. (2009).
\newblock On {Statistical} {Properties} of {Hyperbolic} {Systems}
  with {Singularities}.
\newblock {\em Journal of Statistical Physics}, 136(4):615--642.

\bibitem[Demers and Wojtkowski, 2009]{demers_family_2009}
Demers, M.~F. and Wojtkowski, M.~P. (2009).
\newblock A family of pseudo-{Anosov} maps.
\newblock {\em Nonlinearity}, 22(7):1743--1760.

\bibitem[Franjione et~al., 1992]{franjione_symmetry_1992}
Franjione, J.~G., Ottino, J.~M., and Smith, F.~T. (1992).
\newblock Symmetry concepts for the geometric analysis of mixing flows.
\newblock {\em Philosophical Transactions of the Royal Society of London.
  Series A: Physical and Engineering Sciences}, 338(1650):301--323.

\bibitem[Hertzsch et~al., 2007]{hertzsch_dna_2007}
Hertzsch, J.-M., Sturman, R., and Wiggins, S. (2007).
\newblock {DNA} {Microarrays}: {Design} {Principles} for {Maximizing}
  {Ergodic}, {Chaotic} {Mixing}.
\newblock {\em Small}, 3(2):202--218.

\bibitem[Katok and Strelcyn, 1986]{katok_invariant_1986}
Katok, A. and Strelcyn, J.-M. (1986).
\newblock {\em Invariant {Manifolds}, {Entropy} and {Billiards}. {Smooth}
  {Maps} with {Singularities}}.
\newblock Lecture {Notes} in {Mathematics}. Springer-Verlag, Berlin Heidelberg.

\bibitem[MacKay, 2006]{mackay_cerbelli_2006}
MacKay, R. (2006).
\newblock Cerbelli and {Giona}'s {Map} {Is} {Pseudo}-{Anosov} and {Nine}
  {Consequences}.
\newblock {\em Journal of Nonlinear Science}, 16(4):415--434.

\bibitem[Myers~Hill et~al., 2021]{myers_hill_continuous_2021}
Myers~Hill, J., Sturman, R., and Wilson, M. C.~T. (2021).
\newblock A {Continuous} {Family} of {Non}-{Monotonic} {Toral} {Mixing} {Maps}.
\newblock {\em arXiv:2112.07346}.

\bibitem[Oseledets, 1968]{oseledets_multiplicative_1968}
Oseledets, V.~I. (1968).
\newblock A multiplicative ergodic theorem. {Lyapunov} characteristic numbers
  for dynamical systems.
\newblock {\em Transactions of the Moscow Mathematical Society}, 19:197--231.

\bibitem[Ottino, 1989]{ottino_kinematics_1989}
Ottino, J.~M. (1989).
\newblock {\em The {Kinematics} of {Mixing}: {Stretching}, {Chaos}, and
  {Transport}}.
\newblock Cambridge University Press.

\bibitem[Pesin, 1977]{pesin_characteristic_1977}
Pesin, Y.~B. (1977).
\newblock Characteristic lyapunov exponents and smooth ergodic theory.
\newblock {\em Russian Mathematical Surveys}, 32(4):55.

\bibitem[Przytycki, 1983]{przytycki_ergodicity_1983}
Przytycki, F. (1983).
\newblock Ergodicity of toral linked twist mappings.
\newblock {\em Annales scientifiques de l'École Normale Supérieure},
  16(3):345--354.

\bibitem[Springham and Sturman, 2014]{springham_polynomial_2014}
Springham, J. and Sturman, R. (2014).
\newblock Polynomial decay of correlations in linked-twist maps.
\newblock {\em Ergodic Theory and Dynamical Systems}, 34(5):1724--1746.

\bibitem[Sturman et~al., 2006]{sturman_ottino_wiggins_2006}
Sturman, R., Ottino, J.~M., and Wiggins, S. (2006).
\newblock {\em The Mathematical Foundations of Mixing: The Linked Twist Map as
  a Paradigm in Applications: Micro to Macro, Fluids to Solids}.
\newblock Cambridge Monographs on Applied and Computational Mathematics.
  Cambridge University Press.

\bibitem[Sturman and Springham, 2013]{sturman_rate_2013}
Sturman, R. and Springham, J. (2013).
\newblock Rate of chaotic mixing and boundary behavior.
\newblock {\em Physical review. E, Statistical, nonlinear, and soft matter
  physics}, 87:012906.

\bibitem[Viana, 2014]{viana_lectures_2014}
Viana, M. (2014).
\newblock {\em Lectures on {Lyapunov} {Exponents}}.
\newblock Cambridge {Studies} in {Advanced} {Mathematics}. Cambridge University
  Press, Cambridge.

\bibitem[Wojtkowski, 1980]{wojtkowski_linked_1980}
Wojtkowski, M. (1980).
\newblock Linked {Twist} {Mappings} {Have} the {K}-{Property}.
\newblock {\em Annals of the New York Academy of Sciences}, 357(1):65--76.

\bibitem[Young, 1998]{young_statistical_1998}
Young, L.-S. (1998).
\newblock Statistical {Properties} of {Dynamical} {Systems} with {Some}
  {Hyperbolicity}.
\newblock {\em Annals of Mathematics}, 147(3):585--650.

\bibitem[Young, 1999]{young_recurrence_1999}
Young, L.-S. (1999).
\newblock Recurrence times and rates of mixing.
\newblock {\em Israel Journal of Mathematics}, 110(1):153--188.

\end{thebibliography}
